\documentclass[a4paper]{article}

\usepackage[ruled,norelsize]{algorithm2e}
\usepackage[noend]{algorithmic}
\usepackage{amsmath,amssymb,mathtools,amsfonts,amsthm}
\usepackage{nicefrac}
\usepackage{calc}
\usepackage{subfig}
\usepackage{bbm}
\usepackage{paralist}
\usepackage{enumitem}
\usepackage{color}
\usepackage{graphicx}
\usepackage{booktabs}
\usepackage[hidelinks]{hyperref}
\usepackage[utf8]{inputenc}
\usepackage{scalerel}
\usepackage{fancyhdr}

\newcommand{\N}{\mathbb{N}}
\newcommand{\R}{\mathbb{R}}

\newcommand{\Order}{\mathcal{O}}
\newcommand{\Normal}{\mathcal{N}}

\newcommand{\bmat}{\begin{pmatrix}}
\newcommand{\emat}{\end{pmatrix}}

\title{A Superlinearly Convergent Evolution Strategy}

\fancypagestyle{firststyle}
{
   \fancyhf{}
   \fancyfoot{}
}

\author{
	Tobias Glasmachers\\
	Faculty of Computer Science, Institute for Neural Computation,\\
	Ruhr-University Bochum, Germany\\
	\texttt{tobias.glasmachers@ini.rub.de}
}

\date{}

\begin{document}

\maketitle

\begin{abstract}
We present a hybrid algorithm between an evolution strategy and a quasi
Newton method. The design is based on the Hessian Estimation Evolution
Strategy, which iteratively estimates the inverse square root of the
Hessian matrix of the problem. This is akin to a quasi-Newton method and
corresponding derivative-free trust-region algorithms like NEWUOA. The
proposed method therefore replaces the global recombination step
commonly found in non-elitist evolution strategies with a quasi-Newton
step. Numerical results show superlinear convergence, resulting in
improved performance in particular on smooth convex problems.
\end{abstract}

\section{Introduction}

Evolution strategies (ES)
\cite{rechenberg1973evolutionsstrategie,beyer2001theory} are
evolutionary algorithms specialized in minimizing a black-box function
$f : \R^d \to \R$. They approach the problem by sampling from a
multi-variate normal distribution, which is shifted and scaled online
during the optimization run to match the currently observed problem
characteristics. State-of-the-art variable-metric ES like the Covariance
Matrix Adaptation ES (CMA-ES)
\cite{hansen2001completely,kern2004learning} adjust the full covariance
matrix of the distribution, which allows them to achieve linear
convergence at a problem-independent rate~\cite{gissler2024linear}.

Due to a general result by Teytaud and Gelly \cite{teytaud2006general},
no optimization algorithm based on ranking of solutions can converge
faster than linear. That is the case for most evolution strategies,
which typically converge at a linear rate of $\Theta(1/d)$, where $d$ is
the problem dimension
\cite{jaegerskupper2007algorithmic,akimoto2018drift,morinaga2019generalized,akimoto2022global}.
In other words, the distance to the optimum shrinks by a constant
factor in each iteration.

In contrast, the line of derivative-free optimization algorithms
developed by Powell \cite{powell2006newuoa} can well achieve superlinear
convergence. That becomes possible due to its reliance on function
values, in contrast to being restricted to their ranks. Using ranks has
pros and cons. A clear advantage is invariance to rank-preserving
(strictly monotonically increasing) transformations of objective
function values, and the added robustness that comes with it. The price
to pay is the limitation to a linear convergence rate.

In this paper we propose a hybrid algorithm, called the Quasi-Newton
Evolution Strategy (QN-ES). At its core it is a modern variable metric
evolution strategy, closely related and strongly influenced by CMA-ES
\cite{hansen2001completely,kern2004learning}. It samples offspring from
a multivariate Gaussian distribution and it uses cumulative step size
adaptation for controlling global step lengths. The mechanism for
updating the covariance matrix of the sampling distribution is taken
over from the Hessian Estimation Evolution Strategy (HE-ES)
\cite{glasmachers2020hessian}. In contrast to CMA-ES, that algorithm
benefits from monotonic improvements of its Hessian approximation, which
also results in linear convergence at a problem-independent rate on all
convex quadratic functions~\cite{glasmachers2022convergence}.

CMA-ES and HE-ES use global intermediate recombination for updating
the distribution mean. In our algorithm we replace that mechanism with a
quasi-Newton step of the form
\begin{align}
	p = - B^{-1} g, \label{eq:qn}
\end{align}
where $p$ is the step performed on the distribution mean, $B$ is a
positive definite matrix estimating the Hessian, and $g$ is an estimate
of the gradient.

Mechanisms for obtaining both estimates exist already. Natural Evolution
Strategies (NES) \cite{wierstra2014natural} estimate the gradient of
expected fitness with respect to the distribution parameters. Taking
only the mean vector as a parameter yields an estimate $g$ of the
gradient as a special case. Most NES algorithms introduce rank-based
selection in order to achieve linear convergence, in contrast to
sublinear convergence caused by the na{\"i}ve application of gradient
estimates \cite{beyer2014convergence}. The approach was later formalized
and generalized to what is now known as the Information Geometric
Optimization (IGO) \cite{ollivier2017information} framework.
HE-ES offers an iterative updating mechanism for estimating the Hessian.
While the latter learning mechanism necessarily uses function values (in
contrast to ranks), we revert to the basics of NES to obtain an actual
gradient estimator. A similar route was taken by Krause
\cite{krause2019large} for noisy problems in reinforcement learning.

Our decisive step then is to replace global intermediate recombination
with an update of the form~\eqref{eq:qn}. We show that such an update
can be performed efficiently in the HE-ES framework, without the need to
invert a matrix or solve a system of equations as indicated by
equation~\eqref{eq:qn}.

Replacing the recombination operator in an evolution strategy may seem
to be a minor change. However, its consequences are far-reaching. Taking
an analogy from gradient-based optimization, the difference can be as
significant as between a simple gradient descent and the BFGS algorithm.
In section~\ref{section:experiments} we will present empirical data
indicating that the new update indeed breaks the barrier of linear
convergence, achieving superlinear convergence on smooth convex
problems.

\section{Background}

Before diving into the details of our method we present the necessary
background. We start with the HE-ES algorithm and then turn to
quasi-Newton algorithms.

\subsection{The Hessian Estimation Evolution Strategy}

The Hessian Estimation Evolution Strategy (HE-ES) was first presented in
\cite{glasmachers2020hessian} and thoroughly analyzed in
\cite{glasmachers2022convergence}. Its most important features are:
\begin{compactitem}
\item
	Like CMA-ES, it uses non-elitist selection, keeping (the weighted
	recombination of) $\mu$ out of $\lambda$ offspring for the next
	generation. In other words, the distribution mean $m$ is updated by
	taking a weighted mean of the offspring samples with rank-based
	weights.
\item
	The step size $\sigma$ is adapted with the cumulative stepsize
	adaptation (CSA) rule \cite{ostermeier1994step}. However, other
	rules can easily be incorporated, even success-based methods, since
	HE-ES always evaluates the distribution mean.
\item
	Offspring are drawn from a multivariate normal distribution
	$\Normal(m, \sigma^2 A^T A)$ in the form $x = m + \sigma A b$,
	where $b \sim \Normal(0, I)$ is a standard normal vector. The
	transformation matrix $A$ is a factor of the covariance matrix.
\item
	HE-ES draws dependent Gaussian samples: offspring always come in
	pairs, implementing mirror sampling. Furthermore, the directions
	of pairs are orthogonal, before applying $A$. If the number of
	pairs exceeds the problem dimension $d$ then each batch of $d$
	directions forms an orthogonal basis.
\item
	The update mechanism for the matrix $A$ is the core of the method.
	$A$ is adapted so that $A^T A$ converges to a multiple of the
	inverse of the Hessian of a convex quadratic objective function.
	Care is taken to keep the matrix positive definite, and to avoid
	instability due to unreasonably large changes.
\item
	The mean $m$ encodes the position, the step size $\sigma$ the
	scale, and the transformation $A$ the shape of the distribution.
	The roles are clearly separated by means of the constraint
	$\det(A) = 1$.
\item
	The (1+4)-HE-ES, an elitist variant, was the first evolution
	strategy for which the benefit of covariance matrix adaptation was
	thoroughly proven: \cite{glasmachers2022convergence} established
	linear convergence on convex quadratic problems at a rate that
	depends only on the problem dimension, but not on the problem
	instance.
\end{compactitem}
The pseudo-code of HE-ES is stated in Algorithm~\ref{algorithm:HE-ES}.
The update%
\footnote{
	We corrected a minor inaccuracy in the last line of
	Algorithm~\ref{procedure:computeG}. The original formulation in
	\cite{glasmachers2020hessian} did not result in a matrix with unit
	determinant if the number of pairs exceeds the problem dimension.
	Our correction coincides with the original algorithm for $\tilde
	\lambda \leq d$, which is fulfilled when using the default
	population size. That change has a minimal effect on performance, at
	best. Note that the $\exp(\dots)$ in that line is a matrix
	exponential, which can be computed for example by means of an eigen
	decomposition.}
of the matrix $A$ is computed by Algorithm~\ref{procedure:computeG}.
The following notation from the algorithms will be used later on:
$\tilde \lambda$ is the number of directions or sample pairs, $n_b$ is
the number of batches, $b_{ij}$ is the $i$-th direction of the $j$-th
batch (a standard normal random vector), and $x_{ij}^-$ and $x_{ij}^+$
are the corresponding mirrored samples. The variables $q_{ij}$ are
proportional to the logarithm of the curvature along $b_{ij}$.

\begin{algorithm}
\caption{computeG}
\label{procedure:computeG}
\begin{algorithmic}[1]
\setlength{\itemindent}{0cm}
\STATE{\textbf{input} $b_{ij}$, $f(m)$, $f(x_{ij}^\pm)$, $\sigma$}
\STATE{\textbf{parameters} $\kappa$, $\eta_A$}
\STATE{$h_{ij} \leftarrow \frac{f(x_{ij}^+) + f(x_{ij}^-) - 2 f(m)}{\sigma^2 \cdot \|b_{ij}\|^2}$ \hfill \# estimate curvature along $b_{ij}$}\hspace{2em}
\STATE{\textbf{if} $\max(\{h_{ij}\}) \leq 0$ \textbf{then} \textbf{return} $I$}
\STATE{$c \leftarrow \max(\{h_{ij}\}) / \kappa$}
\STATE{$h_{ij} \leftarrow \max(h_{ij}, c)$ \hfill \# truncate huge steps}\hspace{2em}
\STATE{$q_{ij} \leftarrow \log(h_{ij})$}
\STATE{$q_{ij} \leftarrow q_{ij} - \frac{1}{\tilde \lambda} \cdot \sum_{ij} q_{ij}$ \hfill \# subtract mean $\to$ ensure unit determinant}\hspace{2em}
\STATE{$q_{ij} \leftarrow q_{ij} \cdot \frac{-\eta_A}{2}$ \hfill \# learning rate and inverse square root (exponent $-1/2$)}\hspace{2em}
\STATE{$q_{n_b,j} \leftarrow 0 \quad \forall j \in \{d n_b - \tilde \lambda, \dots, d\}$ \hfill \# neutral update in the unused directions}\hspace{2em}
\STATE{\textbf{return} $\exp \left( \frac{1}{n_b} \sum_{ij} \frac{q_{ij}}{\|b_{ij}\|^2} \cdot b_{ij} b_{ij}^T \right)$}
\end{algorithmic}
\end{algorithm}

\begin{algorithm}
\caption{Hessian Estimation Evolution Strategy (HE-ES)}
\label{algorithm:HE-ES}
\begin{algorithmic}[1]
\STATE{\textbf{input} $m^{(0)} \in \R^d$, $\sigma^{(0)} > 0$, $A^{(0)} \in \R^{d \times d}$}
\STATE{\textbf{parameters} $\tilde \lambda \in \N$, $c_s$, $d_s$, $w \in \mathbb{R}^{2 \tilde \lambda }$}
\STATE{$n_b \leftarrow \lceil \tilde \lambda / d \rceil$}
\STATE{$p_s^{(0)} \leftarrow 0 \in \R^d$}
\STATE{$g_s^{(0)} \leftarrow 0$}
\STATE{$t \leftarrow 0$}
\REPEAT
	\FOR{$j \in \{1, \dots, n_b\}$}
		\STATE{$b_{1j}, \dots, b_{dj} \leftarrow$ \texttt{sampleOrthogonal}()}
	\ENDFOR
	\STATE{$x_{ij}^- \leftarrow m^{(t)} - \sigma^{(t)} \cdot A^{(t)} b_{ij}$ ~~~~~for $i+(j-1)n_b \leq \tilde \lambda$}
	\STATE{$x_{ij}^+ \leftarrow m^{(t)} + \sigma^{(t)} \cdot A^{(t)} b_{ij}$ ~~~~~for $i+(j-1)n_b \leq \tilde \lambda$ \hfill \# mirrored sampling}\hspace{2em}
	\STATE{$A^{(t+1)} \leftarrow A^{(t)} \cdot$ \texttt{computeG}($\{b_{ij}\}$, $f(m)$, $\{f(x_{ij}^\pm)\}$, $\sigma^{(t)}$)} \hfill \# matrix adaptation\hspace{2em}
    \STATE{$w_{ij}^\pm  \leftarrow w_{\text{rank}(f(x_{ij}^\pm))}$}
	\STATE{$m^{(t+1)} \leftarrow \sum_{ij} w_{ij}^\pm \cdot x_{ij}^\pm$} \hfill \# mean update\hspace{2em}
	\STATE{$g_s^{(t+1)} \leftarrow (1 - c_s)^2 \cdot g_s^{(t)} + c_s \cdot (2 - c_s)$}
	\STATE{$p_s^{(t+1)} \leftarrow (1 - c_s) \cdot p_s^{(t)} + \sqrt{c_s \cdot (2 - c_s) \cdot \mu_\text{eff}^\text{mirrored}} \cdot \sum_{ij} (w_{ij}^+ - w_{ij}^-) \cdot b_{ij}$}
	\STATE{$\sigma^{(t+1)} \leftarrow \sigma^{(t)} \cdot \exp\left( \frac{c_s}{d_s} \cdot \frac{\|p_s^{(t+1)}\|}{\chi_d} - \sqrt{g_s^{(t+1)}} \right)$ \hfill \# CSA}\hspace{2em}
	\STATE{$t \leftarrow t + 1$}
\UNTIL{ \textit{stopping criterion is met} }
\end{algorithmic}
\end{algorithm}

There is a simple intuition for the effect of the update matrix $G$: The
algorithm measures the curvature (estimates the second derivative) of
the objective function along the orthogonal set of sampled directions by
assessing the distribution mean and a pair of mirrored samples. At the
same time, the matrix $A$ gives rise to the corresponding curvature of
the quadratic model
$\frac12 x^T (A^T A)^{-1} x = \frac12 \|A^{-1} x\|^2$.
If the quotient of two curvatures deviates between measurement and
quadratic model then a correction is applied that moves the model
curvatures closer to the measurements, where a learning rate of $\eta_A
= 1$ means that the old values are completely overwritten by the
measurements. Aside from this update mechanism and the special
requirements on the samples that come with it, HE-ES shares most
algorithmic components with CMA-ES. For further details we refer the
interested reader to
\cite{glasmachers2020hessian,glasmachers2022convergence}.

\subsection{Quasi-Newton Methods}

A quasi-Newton method is an iterative gradient-based optimizer
performing steps of the form $x_t = x_{t-1} + \alpha p$ with step size
$\alpha$ preferably equal to one and direction
$p = -B^{-1} \nabla f(x_{t-1})$, which is analogous to
equation~\eqref{eq:qn}. $B$ is an approximation of the Hessian
$\nabla^2 f(x_{t-1})$. The corresponding quantities exist if the
objective function is twice continuously differentiable. Quasi-Newton
methods require access to first derivatives, while the Newton method
also requires access to second derivatives.

If $B$ equals the Hessian then the step is a Newton step. It is
meaningful for optimization only if $B$ is positive definite. Newton
steps solve a convex quadratic problem in a single step and otherwise
they enjoy quadratic (very fast) local convergence, but the method can
easily diverge when started too far away from the optimum. Quasi-Newton
methods like the BFGS algorithm are a bit slower, but still converge
superlinearly \cite{nocedal2006numerical}. The intuition is that the
approach based on an approximate quadratic model can come arbitrarily
close to an actual Newton step because the quadratic model converges to
the second order Taylor approximation of the objective function while
the algorithm approaches the optimum.

Several variants of quasi-Newton methods exist. They differ in how they
update the matrix $B$ based on each gradient evaluation. A minimal
requirement and a guiding design principle is that the update must
fulfill the \emph{secant equation}
$$ B (x_t-x_{t-1}) = \nabla f(x_t) - \nabla f(x_{t-1}), $$
which demands that $B$ shall be updated such that after the update the
quadratic model
$$ \frac12 (x-x_t)^T B (x-x_t) + \nabla f(x_t)^T (x-x_t) + f(x_t) $$
predicts the just evaluated gradient correctly. Quasi-Newton methods
perform this update in such a way that it implements a minimal change of
$B$, measured by a suitable matrix norm. Important instances of
quasi-Newton algorithms are the DFP, SR1, and BFGS update rules
\cite{nocedal2006numerical}. Of these, all but the SR1 rule keep the
matrix $B$ positive definite, even if the underlying Hessians
$\nabla^2 f(x_t)$ are indefinite.

Quasi-Newton methods can be combined with a trust region approach. Such
algorithms are very successful in the derivative free optimization (DFO)
domain, in particular the {NEWUOA} algorithm \cite{powell2006newuoa} and
its successors. In contrast to gradient-based algorithms, these methods
need to carefully balance two objectives, namely fast progress by means
of quasi-Newton steps, and a well-spread sample set for avoiding
numerical collapse due to ill-conditioning. A detailed discussion of
these techniques is found in~\cite{conn2009introduction}.

\section{The Quasi Newton Evolution Strategy}

Now we have all ingredients in place to present our novel Quasi-Newton
Evolution Strategy (QN-ES). Compared with the pseudo code of
Algorithm~\ref{algorithm:HE-ES}, the only change is the mean update in
line~14. Furthermore, we demand $\tilde \lambda \geq d$, or put
differently, that the directions $b_{ij}$ span the whole space. In fact,
we always make $\tilde \lambda$ an integer multiple of~$d$.

\subsection{Gradient and Quasi-Newton Step}

Quasi-Newton methods build a (local) quadratic model of the objective
function. While the matrix $A$ gives rise to the quadratic term
$\frac12 x^T (A^T A)^{-1} x$ (but see below for a correction), the
gradient estimator $g$ yields the linear term $g^T x$. It is understood
that $g$ shall be estimated from the offspring and their function
values, as these are readily available.

Using the distribution of the samples, it is easy to construct an
estimator of the gradient $\nabla f(A x)$ from finite differences of
mirrored samples:
$$ \nabla f(Ax) \approx \frac{1}{2 \sigma \cdot n_b} \sum_{j=1}^{n_B} \sum_{i=1}^{d} \left(f(x_{ij}^+) - f(x_{ij}^-)\right) \frac{b_{ij}}{\|b_{ij}\|} =: \delta $$
Due to the mirrored sampling approach this is a \emph{central
difference} estimator, which can be very accurate. The construction
obviously requires $\tilde \lambda \geq d$, otherwise we would estimate
a randomly projected gradient.

We can turn $\delta$ into an estimator of the gradient of $f$ by means of
$\nabla f(x) \approx A^{-1} \delta$. The resulting quasi-Newton step is
\begin{align}
	p = -(\eta A^T A) A^{-1} \delta = -\eta A^T \delta,   \label{eq:qn2}
\end{align}
with $B^{-1} = \eta A^T A$. Here, $\eta$ represents the ``global''
curvature. Recall that HE-ES only estimates \emph{relative} curvature
differences between different directions, not the absolute scale. It
maintains the invariant $\det(A) = 1$. That corresponds exactly to
learning the shape of the sampling distribution, in contrast to its
global scale, which is controlled by $\sigma$. However, $\sigma$ plays a
different role: it controls the distance of the samples from the mean,
not the global curvature.

From $\det(A^T A) = 1 = \det((A^T A)^{-1})$ we see that the factor
$\eta$ is required for representing a general Hessian with
$\det(H) \not= 1$. It happens to be the expectation of the means that
are subtracted in line~8 of Algorithm~\ref{procedure:computeG},
normalized by the (squared) step size. We can therefore estimate $\eta$
by keeping track of the subtracted means of the logarithmic curvatures
$q_{ij} - \log(\sigma^2)$. In our implementation, we use an average over
a fixed window of the last $20$ values.

A quasi-Newton step can shrink the distance to the optimum by several
orders of magnitude. CSA-based step size adaptation cannot keep up with
such large progress. We therefore set $\sigma$ to the minimum of the
value proposed by CSA and $\eta \cdot \|\delta\|$.

\subsection{Switch between Recombination and Quasi-Newton}

In early iterations, the quadratic model and in particular the
transformation matrix $A$ cannot be expected to be accurate estimates.
That is for multiple reasons: the algorithm did not have enough data for
estimating the Hessian, and for non-quadratic problems the underlying
Hessian may still change quickly. In that situation, a quasi-Newton step
can be arbitrarily bad, resulting in seemingly erratic steps far away
from the optimum if the initial guess of $B$ underestimates the
curvature in at least one direction. The situation can be addressed by
line search or trust region approaches. Because QN-ES is based on HE-ES,
we already have a well-working alternative update mechanism available.
Therefore, we design an adaptive switching mechanism that performs
either weighted recombination (line~14 in
Algorithm~\ref{algorithm:HE-ES}) or the quasi-Newton
step~\eqref{eq:qn2}. The downside of testing both steps is that one of
the function evaluations gets discarded.

We therefore aim for a switching mechanism that tests both step types
only when necessary. To this end, QN-ES maintains an exponentially
fading record of the rate $R$ at which that the
quasi-Newton step outperforms global weighed recombination. It performs
a recombination step with probability $5/2 \cdot (1-R)$ and a
quasi-Newton step with probability $5/2 \cdot R$, both clipped to the
interval $[\frac{1}{100}, 1]$. Note that at least one of these
probabilities is always one by design, so that a step can always be
taken.

In the cases where both steps are active, the better of the two
(according to the resulting function value) is accepted, and $R$ is
replaced with $0.8 R$ if recombination wins and with $0.8 R + 0.2$ if
the quasi-Newton step wins.

The source code of QN-ES is available from the author's website.
%% TODO: do it!!!

\subsection{Where to place QN-ES as a Hybrid?}

The QN-ES algorithm is designed as a hybrid method. In our
understanding, it is mainly an evolution strategy. However, it can also
be understood as a DFO algorithm in very rough analogy to Powell's
NEWUOA method. In that perspective, two properties are clearly
non-standard: first of all, due to its fallback to a robust mean update,
QN-ES does not require a trust-region approach. Furthermore, the vast
majority of the fitness evaluations is used for robustly estimating a
non-degenerate model, in particular in high dimensions. That can be
inefficient when solving relatively smooth and easy problems. On a
conceptual scale between DFO and ES methods, it is probably fair to
place QN-ES much closer towards the ES side. Yet, the novel approach can
be a first cornerstone of a bridge between the two domains.

\section{Experimental Evaluation}
\label{section:experiments}

% main competitor: HE-ES; secondary competitor: Powell's method
Our assessment of the properties of QN-ES is empirical in nature.
Therefore we present a large battery of experimental results, assessing
strengths and limitations of QN-ES. We put QN-ES into perspective by
comparing its behavior and its performance with two baseline methods at
the opposite ends of the algorithmic spectrum on which QN-ES forms a
hybrid: HE-ES, from which large parts of QN-ES are derived, and Powell's
NEWUOA algorithm, a state-of-the-art DFO method. We rely on the
\texttt{pdfo} package \cite{ragonneau2023pdfo}, which is based on
Powell's original Fortran code. We ran it with default setting with one
exception: the option \texttt{radius\_final} was set to $10^{-100}$ to
avoid early stops caused by hitting a bound on the trust region radius.

\begin{table}
\begin{center}
	\begin{tabular}{c|c}
	\textbf{name} & \textbf{definition} \\
	\hline
	sphere & $f(x) = \|x\|^2$ \\
	ellipsoid & $f(x) = \sum_{i=1}^d 10^{6(i-1)/(d-1)} x_i^2$ \\
	discus & $f(x) = 10^6 x_1^2 + \sum_{i=2}^d x_i^2$ \\
	cigar & $f(x) = x_1^2 + 10^6 \cdot \sum_{i=2}^d x_i^2$ \\
	Rosenbrock & $f(x) = \sum_{i=1}^{d-1} 100 \cdot \left(x_{i+1} - 2 x_i - x_i^2\right)^2 + x_i^2$ \\
	log-sphere & $f(x) = \log(\|x\|^2)$ \\
	one-norm & $f(x) = \sum_{i=1}^d |x_i|$ \\
	sum of different powers~ & $f(x) = \sqrt{\sum_{i=1}^d |x_i|^{2+4(i-1)/(d-1)}}$ \\
	happycat & ~$f(x) = \big((\|x\|^2 - d)^2\big)^{1/4} + \frac1n \left( \frac12 \|x\|^2 + \sum_{i=1}^n x_i \right) + \frac12$ \\
	\end{tabular}\vspace*{0.5em}
	\caption{
		\label{table:benchmarks}
		Definition of the benchmark problems used in this study.
		Note that our formulation of the Rosenbrock function is shifted
		so that its optimum is located at the origin. We use a power of
		$\alpha=1/4$ instead of the originally proposed value of $1/8$
		for the happycat function.
	}
\end{center}
\end{table}

The definitions of the nine benchmark problems used in the experiments
are listed in Table~\ref{table:benchmarks}. Furthermore, we ran
large-scale experiments on the BBOB benchmark
suite~\cite{hansen2010comparing}.

\subsection{Convergence Speed}
% ideal case: sphere
% convex quadratic: ellipsoid, discus, cigar; solve these (and Rosenbrock) to very high precision

\begin{figure}
\begin{center}
	\includegraphics[width=0.32\textwidth]{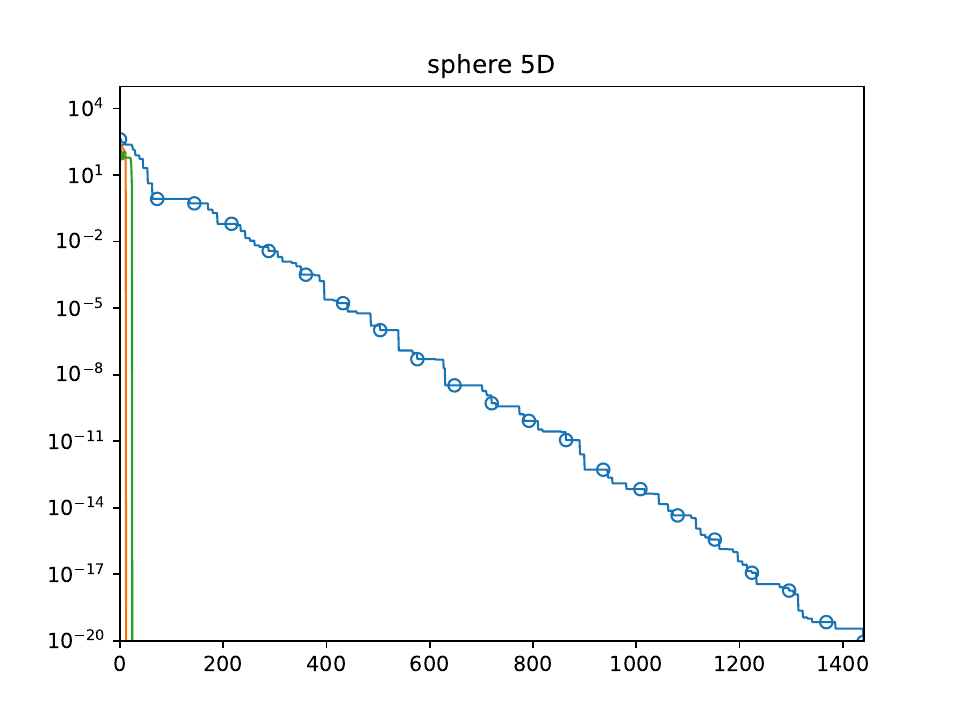}%
	\includegraphics[width=0.32\textwidth]{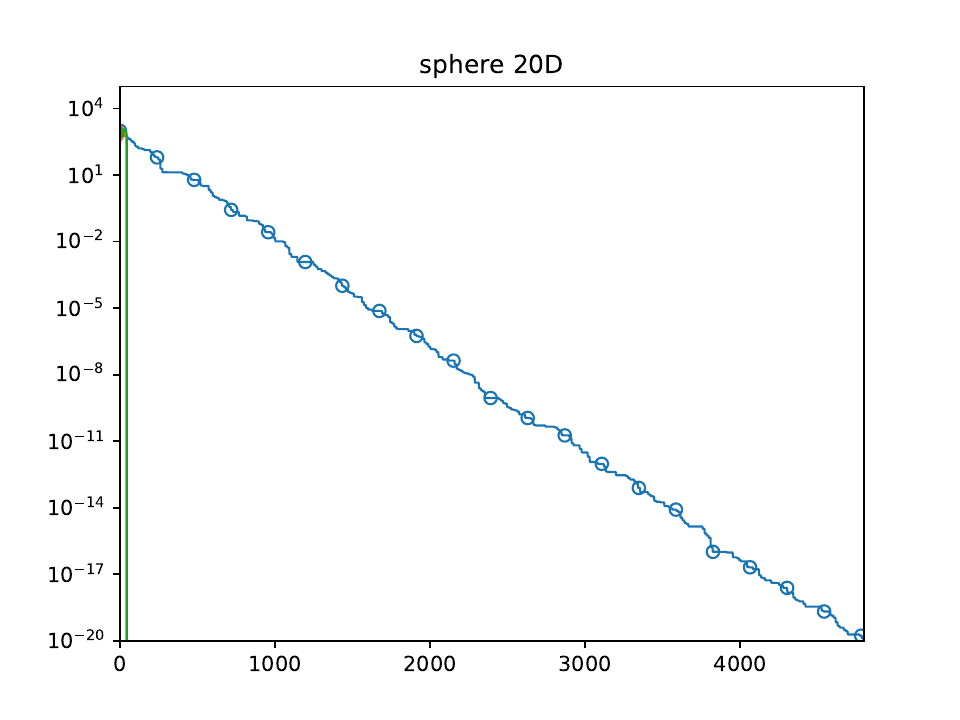}%
	\\
	\includegraphics[width=0.32\textwidth]{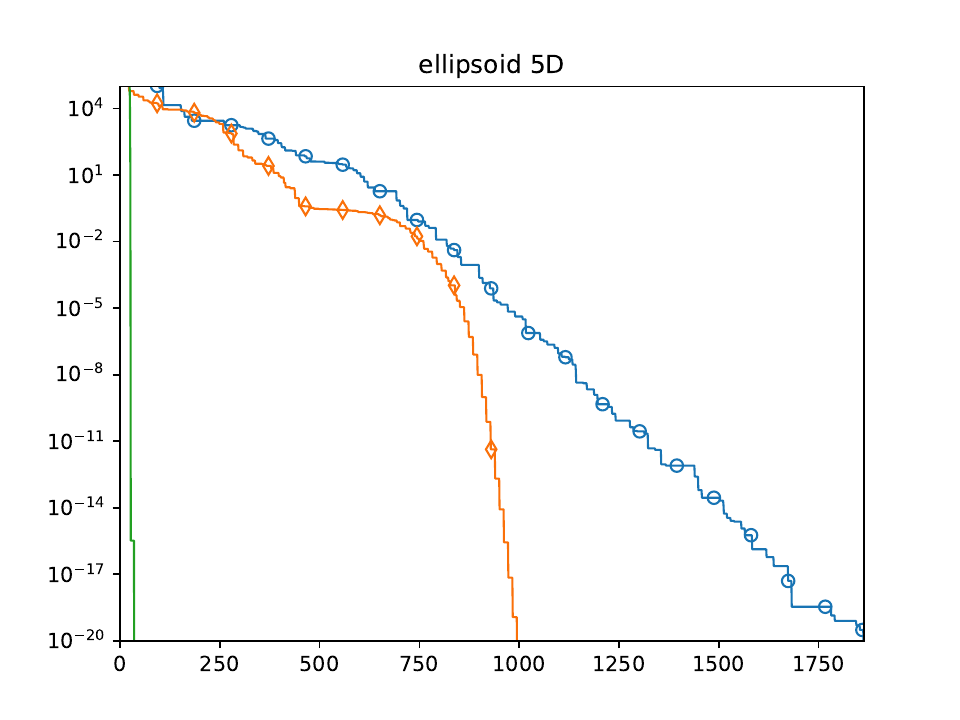}%
	\includegraphics[width=0.32\textwidth]{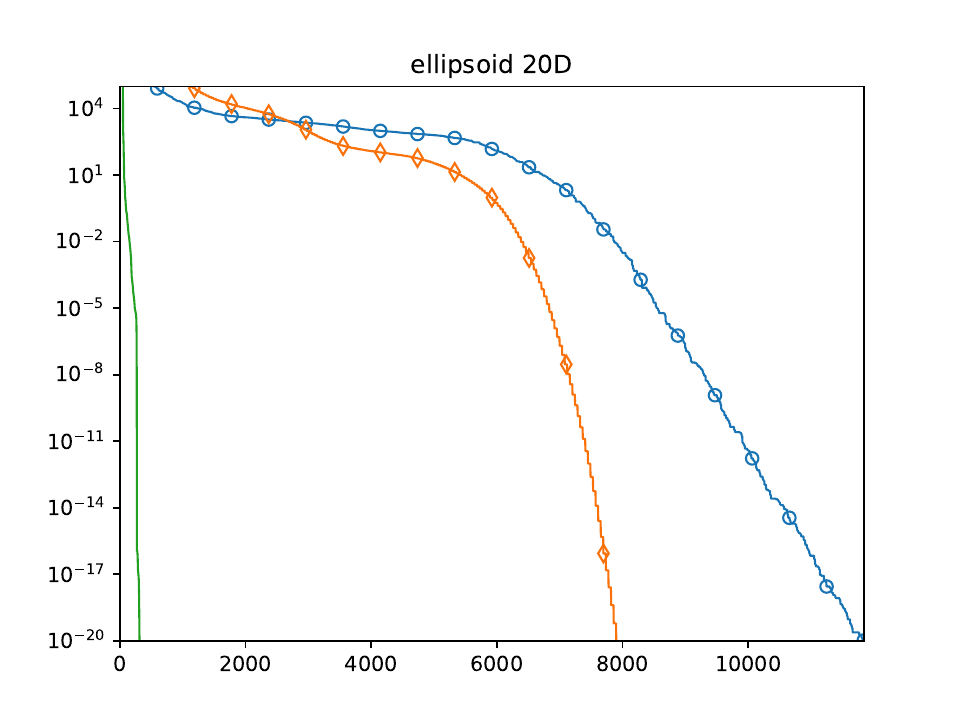}%
	\\
	\includegraphics[width=0.32\textwidth]{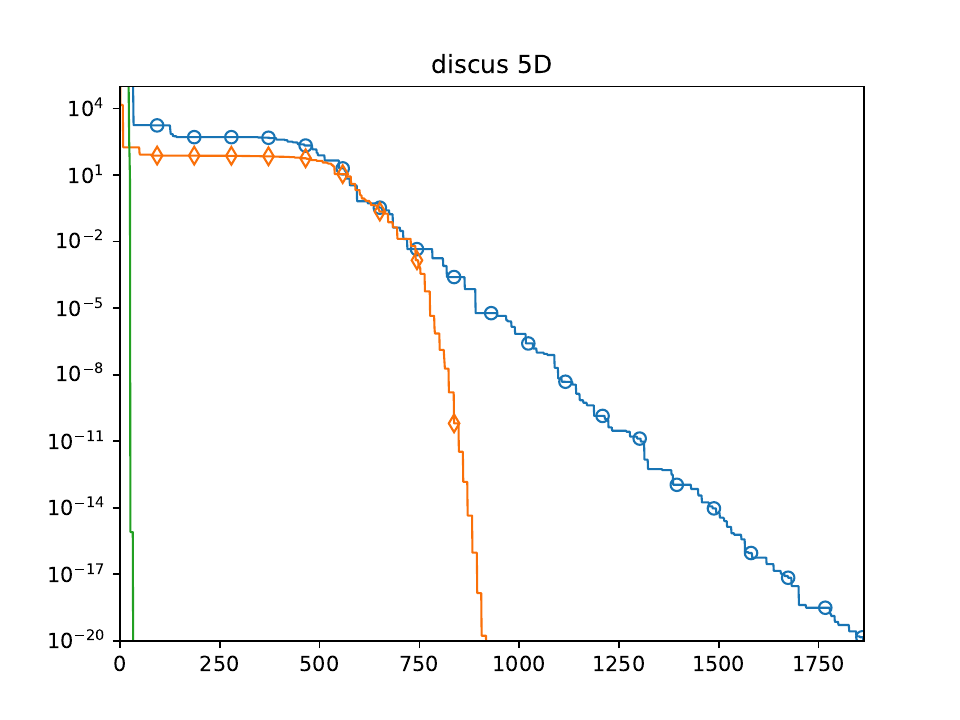}%
	\includegraphics[width=0.32\textwidth]{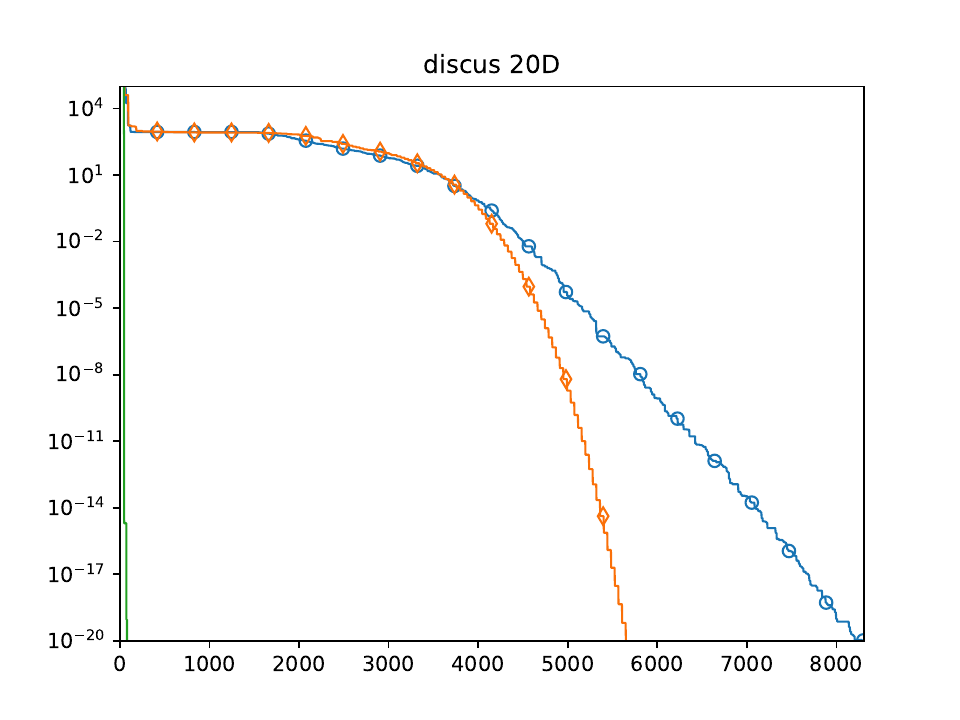}%
	\\
	\includegraphics[width=0.32\textwidth]{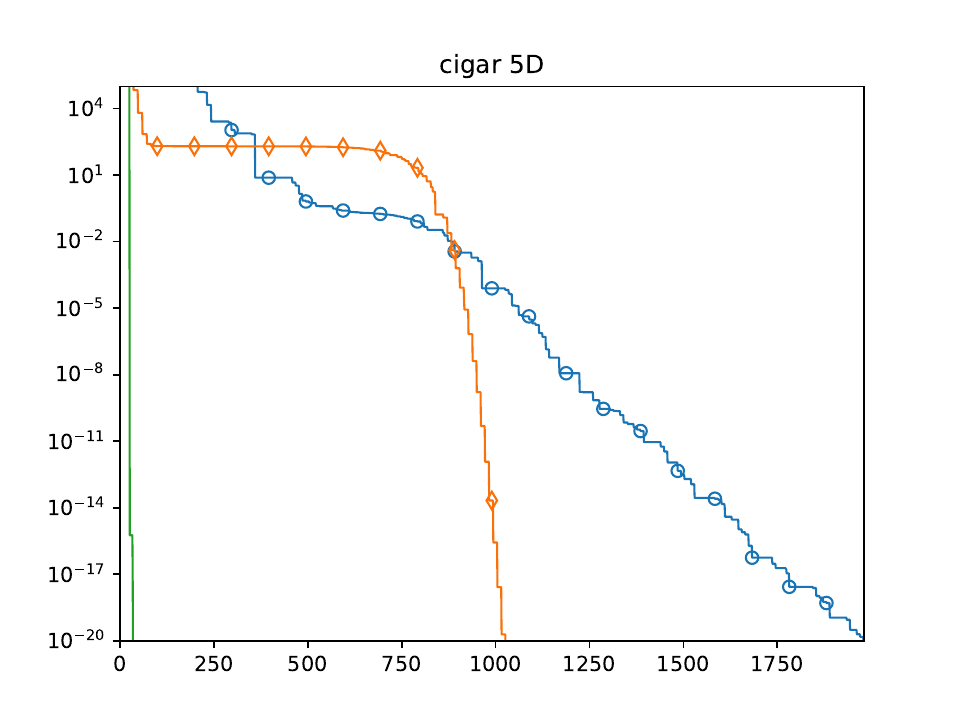}%
	\includegraphics[width=0.32\textwidth]{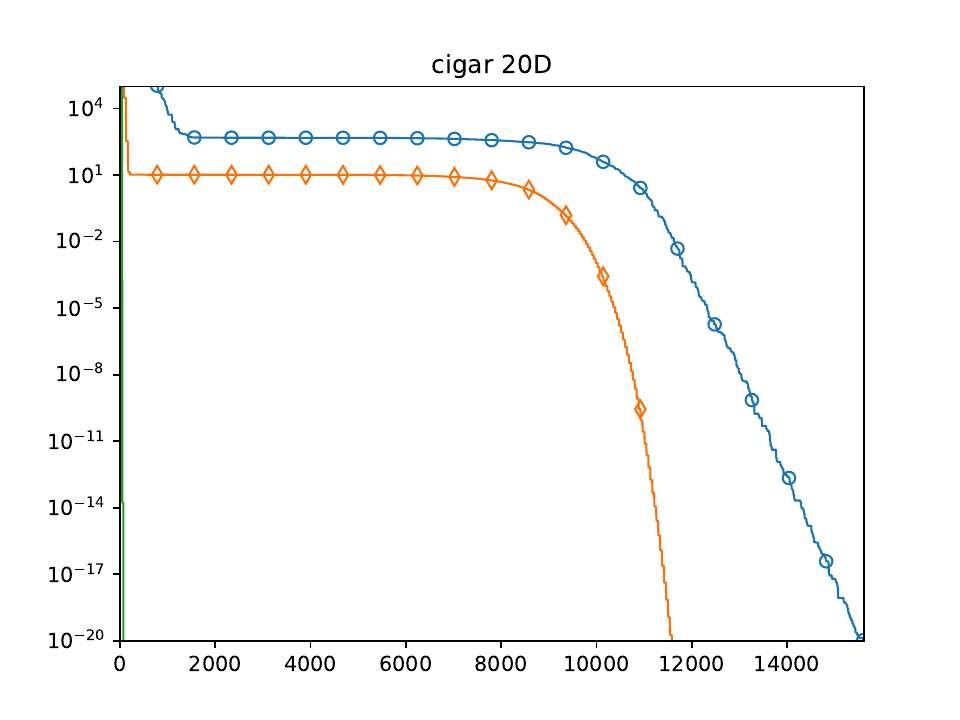}%
	\\
	\includegraphics[width=0.32\textwidth]{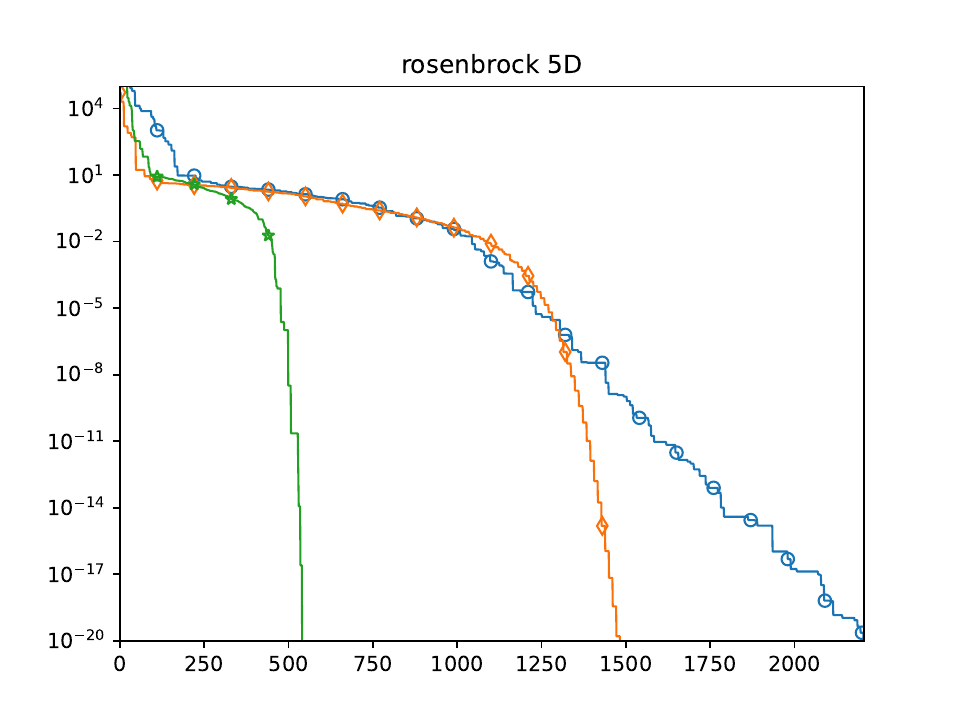}%
	\includegraphics[width=0.32\textwidth]{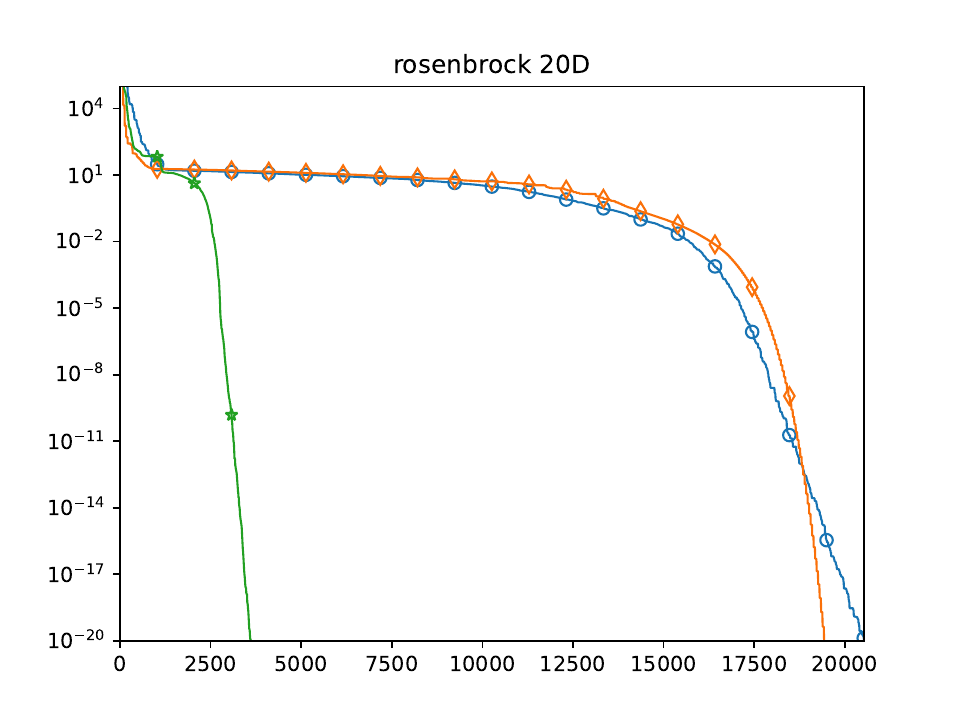}%
	\caption{
		\label{figure:easy}
		Convergence plots (best-so-far $f(x) - f(x^*)$ over the number of
		function evaluations) of Powell's method (green stars), HE-ES (blue
		circles), and QN-ES (orange diamonds) for the problems sphere,
		ellipsoid, discus, cigar, and Rosenbrock (top to bottom) in
		dimensions 5 (left) and 20 (right). The plots show typical single
		runs. In order to assess convergence behavior the problems are
		solved to a high precision of $f(x) - f(x^*) \leq 10^{-20}$.
	}
\end{center}
\end{figure}

We start out by assessing the convergence speed under ideal conditions.
To this end, we consider the simplest non-trivial objective function,
the sphere function. We then move on to other convex quadratic functions
with a high condition number. Finally, we investigate the Rosenbrock
function, a non-convex fourth-order polynomial. Figure~\ref{figure:easy}
shows the corresponding convergence plots.

\begin{figure}
\begin{center}
	\includegraphics[width=0.6\textwidth]{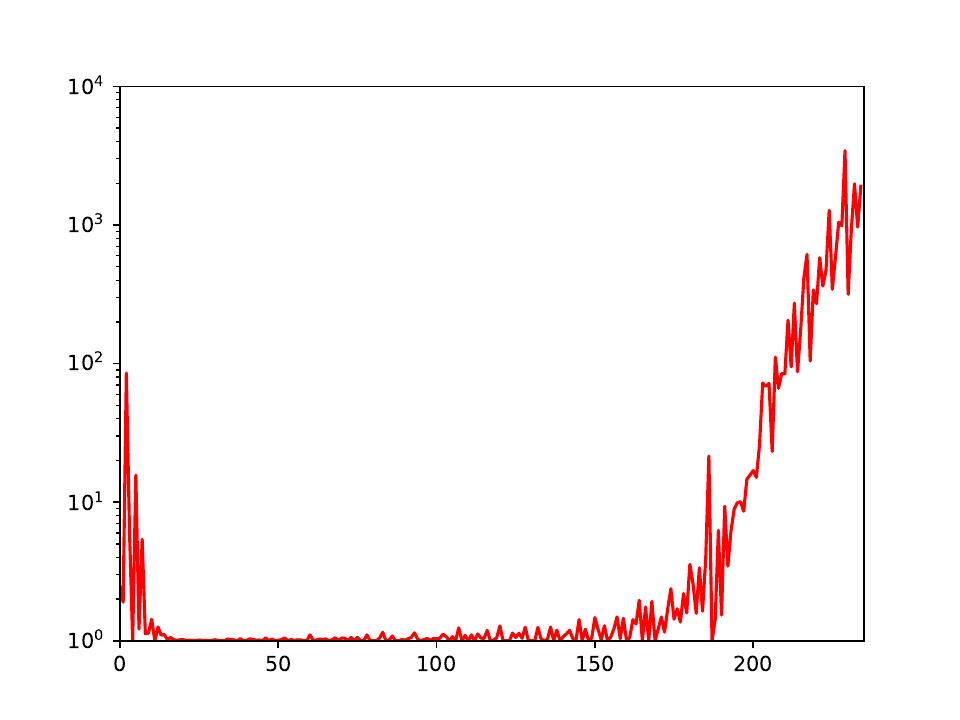}
	\caption{
		\label{figure:progress}
		Multiplicative progress of the function value per iteration of
		QN-ES when minimizing a 10-dimensional Rosenbrock function to
		very high precision. Note the logarithmic scale of the vertical
		axis. After iteration 170, superlinear convergence (a steady
		increase of the convergence rate) is clearly visible. For this
		experiment, the optimum of the Rosenbrock function was shifted
		to the origin in order to make optimal use of floating point
		precision.
	}
\end{center}
\end{figure}

Convex quadratic problems are the best case for all algorithms under
study, but Powell's method and QN-ES profit most due to the perfect
match of (quasi)-Newton steps to the problem class. It is therefore not
surprising that HE-ES performs worst. HE-ES is outperformed by QN-ES on
all quadratic problems. We also observe that the excessive sampling
strategy of using $2 \cdot d$ function evaluations for model building
and only one evaluations per iteration for a quasi-Newton step is
wasteful on these rather simple problems, as demonstrated by the huge
gap between Powell's method and QN-ES. Unsurprisingly, Powell's method
solves quadratic problems extremely quickly to high precision. It
usually takes a number of function evaluations that is only a small
multiple of the problem dimension. QN-ES shows its advantage over HE-ES
in the late convergence phase after successfully learning the
transformation matrix~$A$. This advantage is more pronounced for
low-dimensional problems. The reason is that learning $A$ requires
$\Omega(d)$ function evaluations (rather $\Order(d^2)$ in the worst
case, since the matrix has $d(d+1)/2$ degrees of freedom), while solving
the problem to fixed precision scales only linear in the problem
dimension~$d$.

Powell's method and QN-ES solve the sphere problem nearly instantly to
high precision. On the ill-conditioned problem, QN-ES takes essentially
the same time as HE-ES to learn the transformation $A$, but thereafter
the convergence speed of QN-ES resembles that of Powell's method more
than that of HE-ES. Since we are interested in convergence speed, we
investigate the behavior in more detail. To this end, we analyze the
evolution of $f(x) - f(x^*)$, and in particular its multiplicative
change per iteration. Linear convergence means that the quantity decays
by a constant factor (the convergence rate) per iteration. The results
of a run on the Rosenbrock function is displayed in
Figure~\ref{figure:progress}. It clearly shows that the progress factor
grows to quite large values. A decrease by a factor of more than $10^3$
in a single iteration using only $2d+1$ samples is essentially
impossible with recombination. We therefore clearly see the benefit of
quasi-Newton steps. Furthermore, the factor keeps growing, which means
that the convergence is superlinear.

\subsection{Challenges}

\begin{figure}
\begin{center}
	\includegraphics[width=0.4\textwidth]{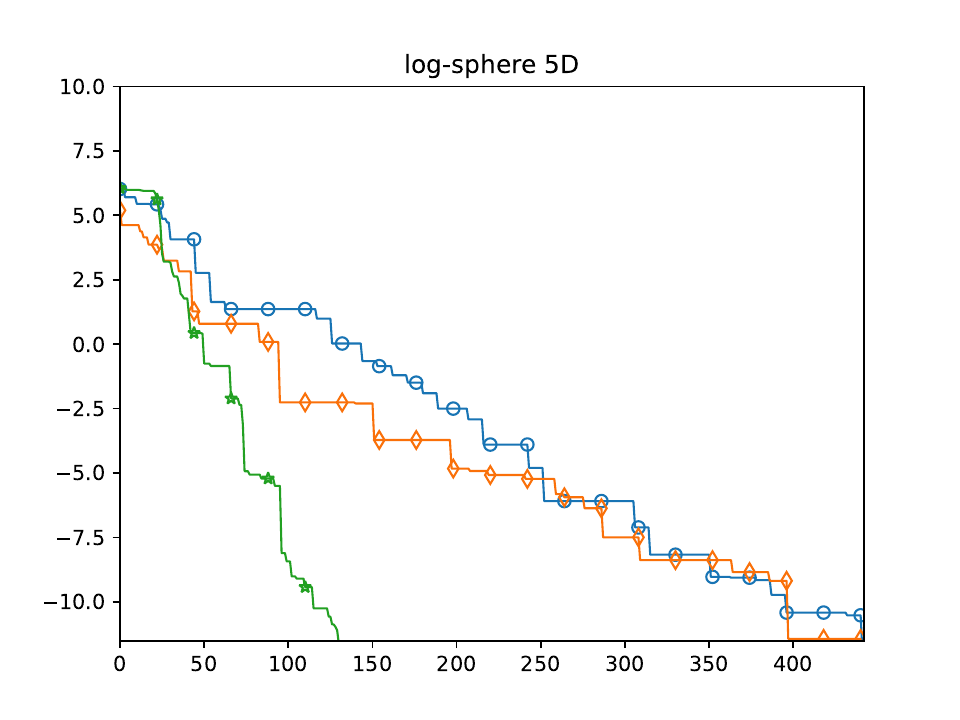}%
	\includegraphics[width=0.4\textwidth]{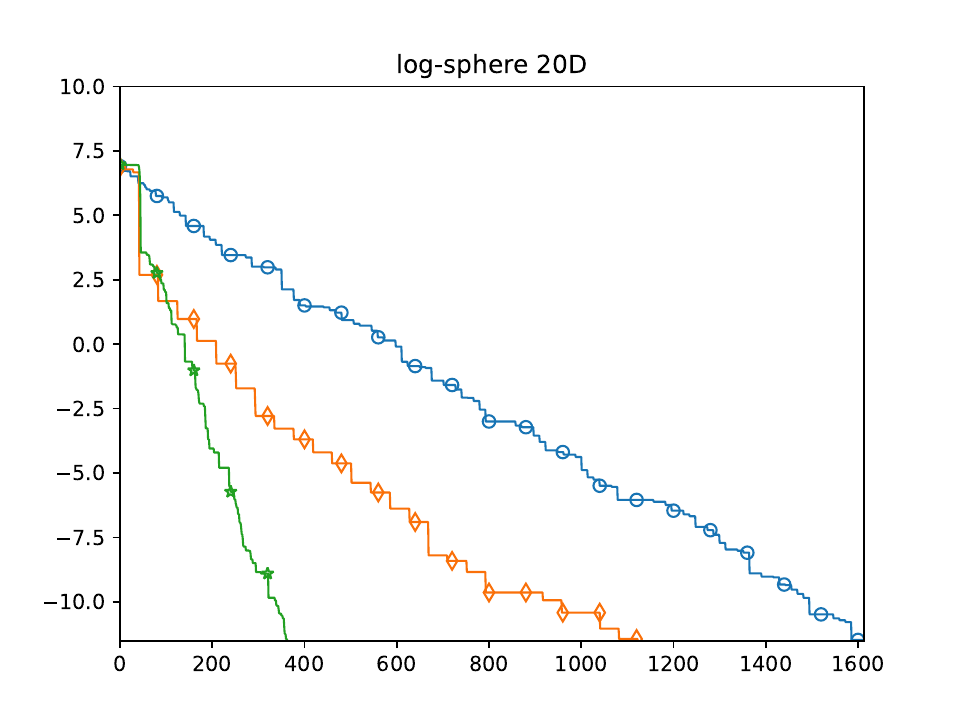}%
	\\
	\includegraphics[width=0.4\textwidth]{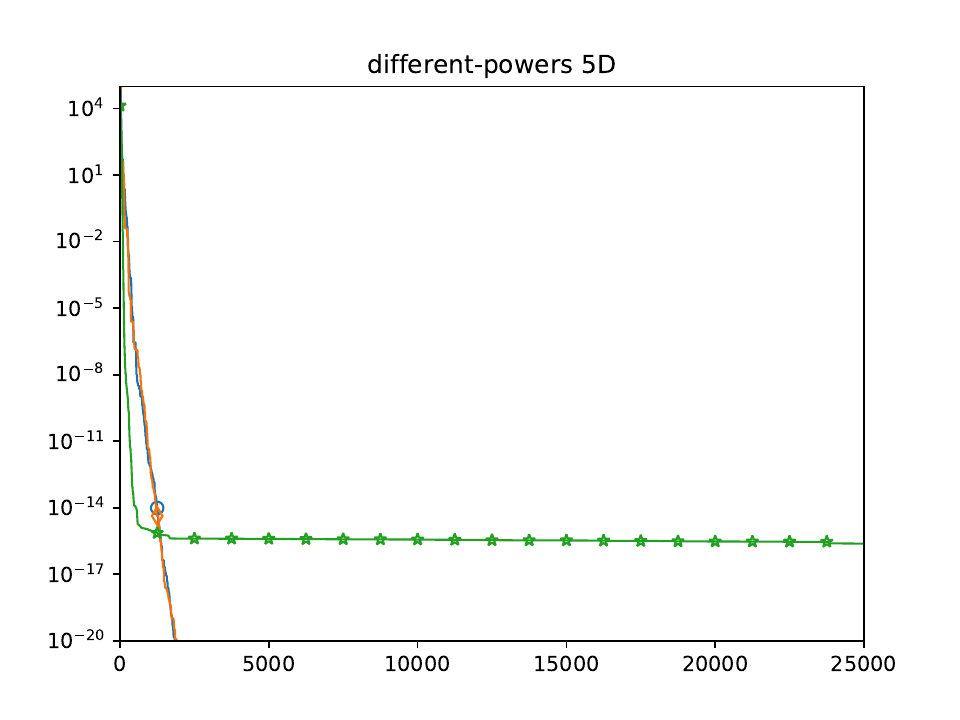}%
	\includegraphics[width=0.4\textwidth]{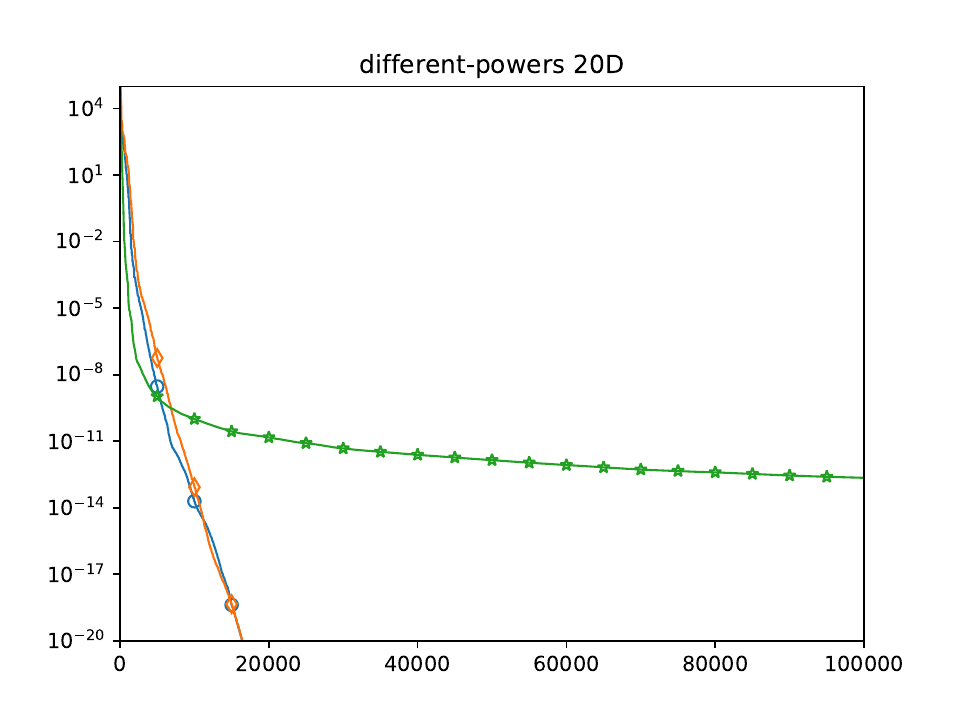}%
	\\
	\includegraphics[width=0.4\textwidth]{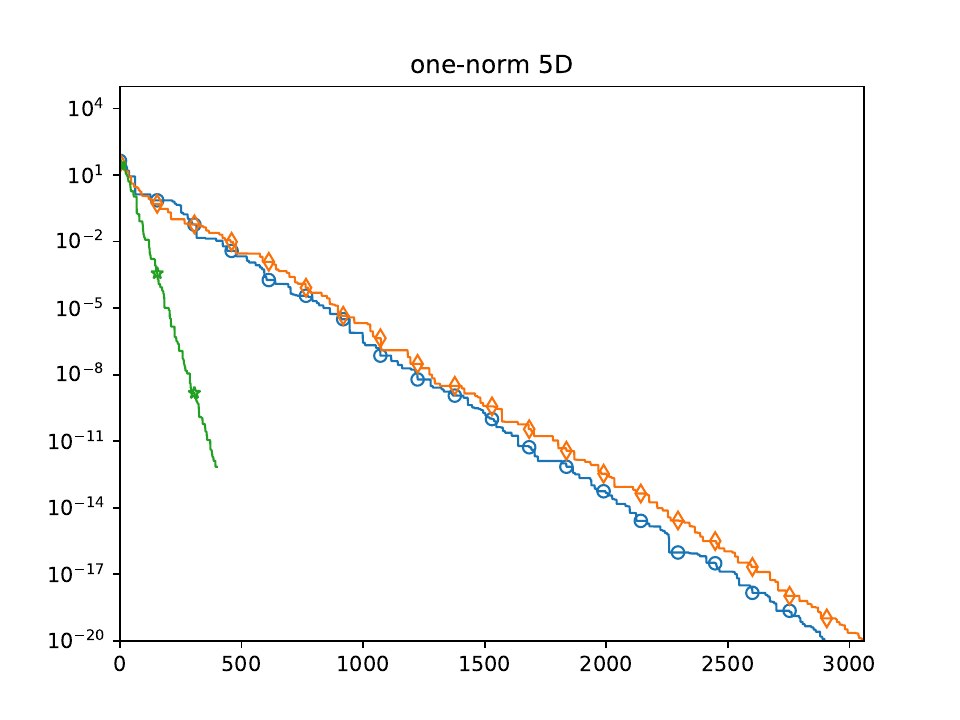}%
	\includegraphics[width=0.4\textwidth]{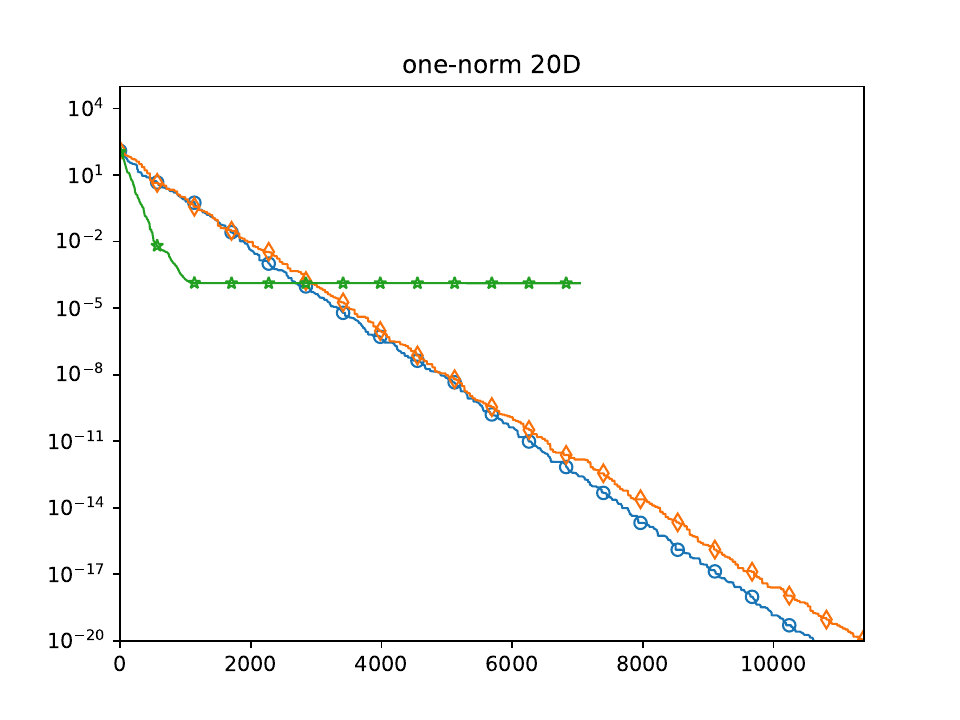}%
	\\
	\includegraphics[width=0.4\textwidth]{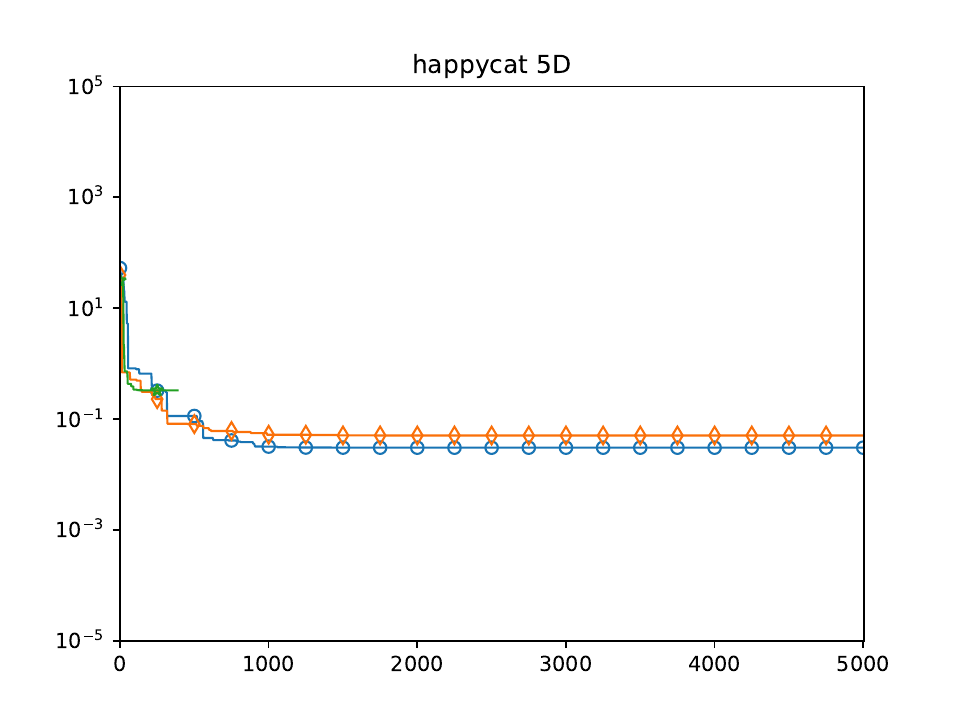}%
	\includegraphics[width=0.4\textwidth]{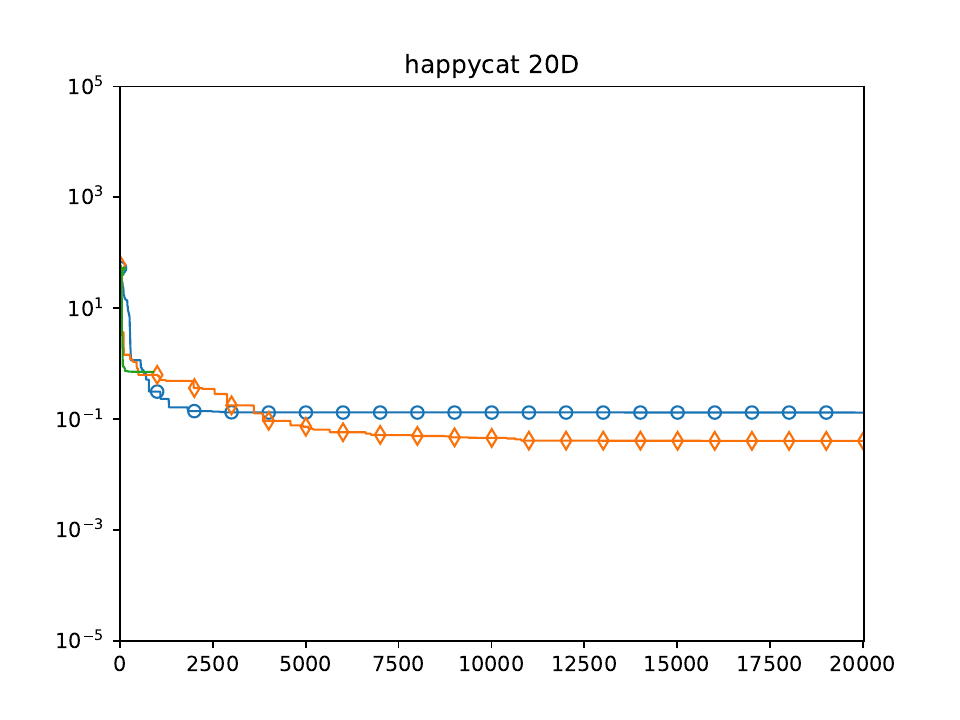}%
	\caption{
		\label{figure:hard}
		Convergence plots of Powell's method (green stars), HE-ES (blue
		dots), and QN-ES (orange diamonds) for the functions log-sphere, sum
		of different powers, one-norm, and happycat (top to bottom) in
		dimensions 5 (left) and 20 (right). The plots show single runs.
	}
\end{center}
\end{figure}

We move on to various challenges, like optimizing non-convex problems,
non-smooth problems with ridges, and extremely ill-conditioned problems.
Most of these are supposedly solved robustly by evolution strategies.
However, since QN-ES is not a typical ES, these experiments are designed
to assess potential limitations of QN-ES. The corresponding results are
found in Figure~\ref{figure:hard}.

The logarithmic sphere function $f(x) = \log(\|x\|^2)$ is an extremely
easy problem for rank-based algorithms because its level set structure
(and therefore the ranking) coincides with that of the sphere function.
However, for QN-ES it is a worst-case scenario because the function is
concave and unbounded in the direction towards the optimum, which makes
Newton steps inapplicable due to the negative curvature. Although HE-ES
is not fully invariant to the logarithmic transformation because it uses
function values for updating $A$, its performance is essentially
unaffected. In contrast, Powell's method as well as QN-ES are heavily
affected by the concave structure, and QN-ES is even outperformed by
HE-ES. In fact, QN-ES systematically falls back to recombination steps
on this problem. We conclude that QN-ES is much more affected by
non-convexity than HE-ES.

The sum-of-different-powers problem is convex. Its difficulty is the
ever increasing conditioning while approaching the optimum. We observe
that Powell's method is initially fast, but it hits a barrier as soon as
the problem conditioning approaches the numerical precision, while both
evolution strategies keep converging. HE-ES and QN-ES show very similar
performance. QN-ES uses quasi-Newton steps only rarely because its
quadratic model keeps dragging behind the changing requirements, and
because the model becomes unreliable after the conditioning of $A$
exceeds the numerical precision. In this case, the switching mechanism
of QN-ES preserves the robustness inherited from HE-ES.

Minimizing a one-norm is a classic problem that appears for example in
machine learning when applying L1-regularization, usually with the goal
of obtaining a sparse solution. The function is continuous and piecewise
linear, forming ridges are the non-smooth transitions. By being
non-differentiable it violates the assumptions of all three algorithms.
QN-ES is (slightly) outperformed by HE-ES because the former relies more
heavily on these assumptions. Powell's method either stalls or fails
completely.

Finally, we included the happycat function \cite{beyer2012happycat},
which was designed as an example of an essentially unsolvable problem
despite being unimodal, continuous, and smooth outside of a circle. It
features a curved and very sharp ridge along which the algorithm must
follow a quadratically decaying trend. It is so delicate that 3D plots
of its graph regularly show artifacts close to the ridge. We set the
power coefficient to $\alpha=1/4$ to make it slightly more regular than
the original version, but the principal difficulty remains. We observe
that all three algorithms fail. Powell's method does not work at all on
this problem. The two evolution strategies keep making progress a bit
longer before stalling far from the optimum.

\subsection{Black-Box Optimization Benchmarking}

We ran QN-ES with a budget of $10^5 \cdot d$ function evaluations per
problem on the Black-Box Optimization Benchmarking (BBOB) suite of the
COmparing Continuous Optimizers (COCO) platform
\cite{hansen2010comparing}. BBOB contains 24 scalable problems with five
instances per problem in dimension 2, 3, 5, 10, 20 and 40. The suite
makes the interpretation of the results convenient by grouping the
problems into meaningful categories.

Several BBOB function categories contain multi-modal problems. We
therefore equipped QN-ES with the same IPOP restart mechanism
\cite{auger2005restart} that was used for HE-ES in
\cite{glasmachers2020hessian}. It restarts the algorithm with doubled
offspring population size (doubling $\tilde \lambda$ and $n_b$) as soon
as objective values stall. We also restart with the same mechanism if
QN-ES aborts due to numerical problems. That wrapper algorithm was
applied to all BBOB problems, not only to the multi-modal ones. For the
competitor methods HE-ES and NEWUOA we downloaded performance data from
the COCO data archive. For a fair comparison, the benchmarked version of
NEWUOA includes a restart mechanism.

\begin{figure}
\begin{center}
	\includegraphics[width=0.333\textwidth]{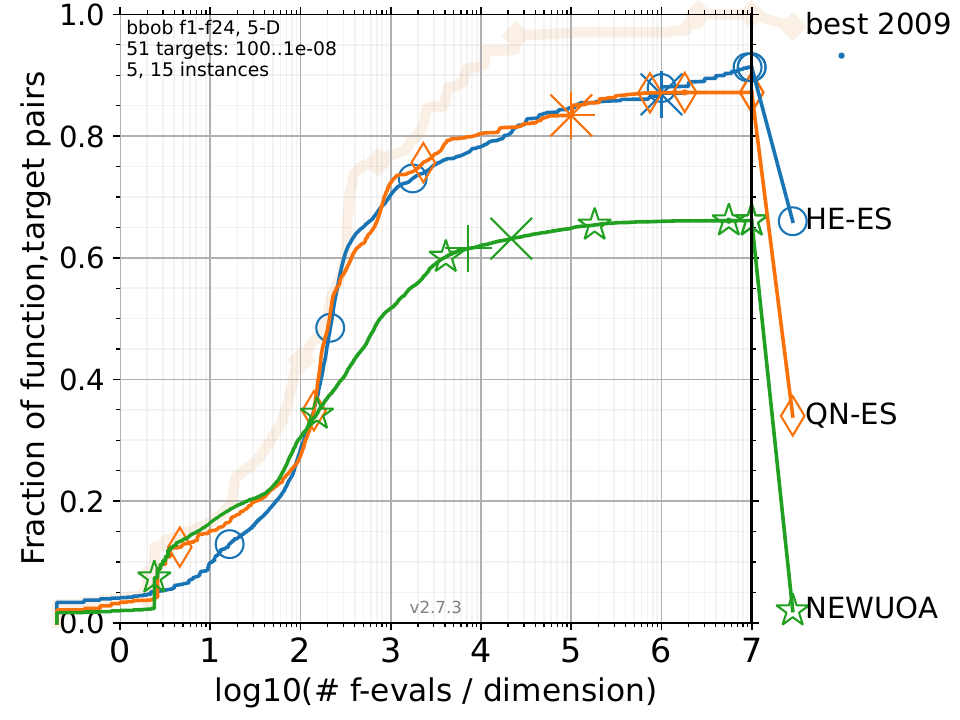}%
	\includegraphics[width=0.333\textwidth]{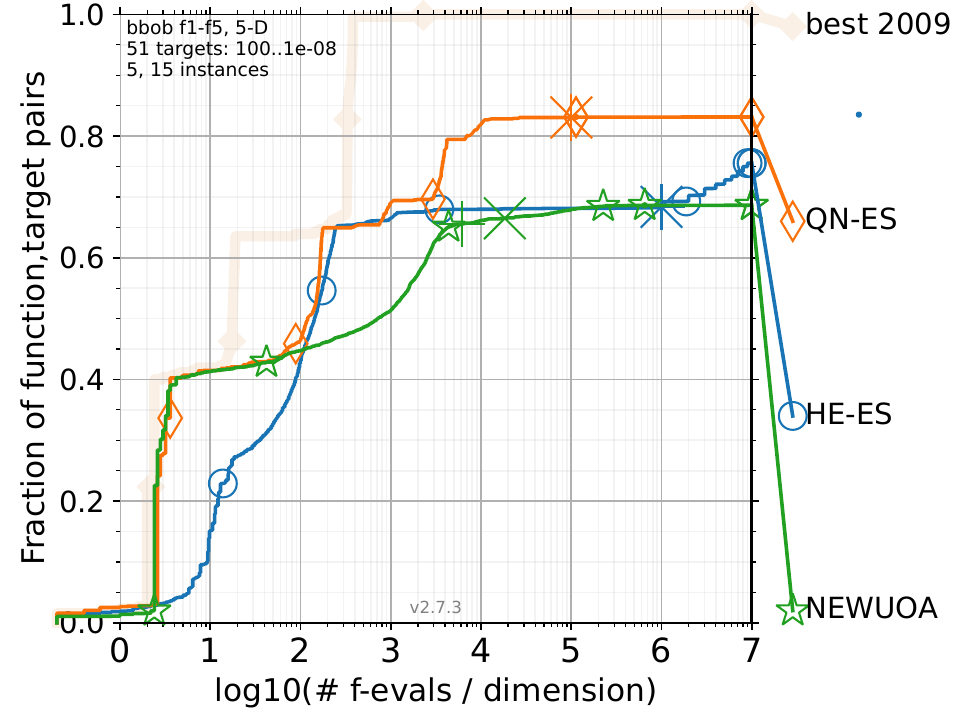}%
	\includegraphics[width=0.333\textwidth]{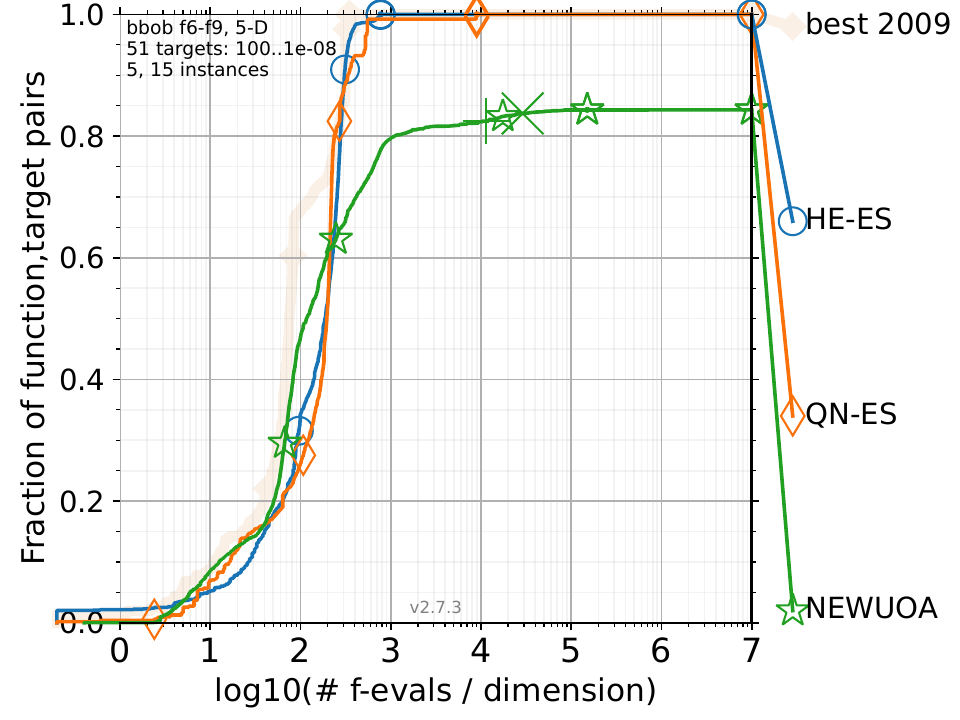}%
	\\
	\includegraphics[width=0.333\textwidth]{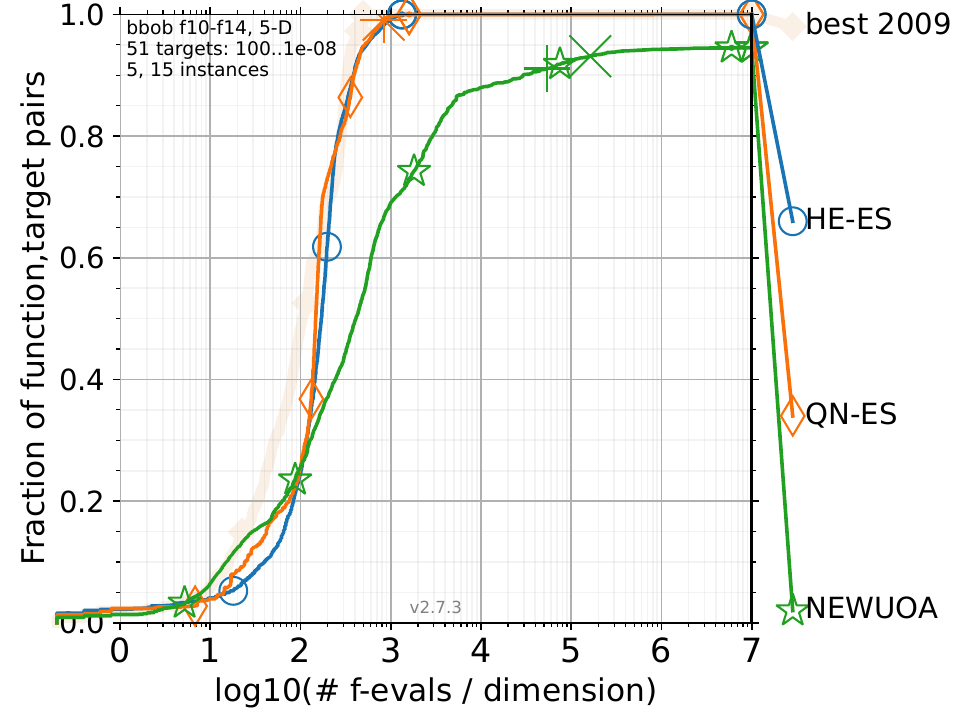}%
	\includegraphics[width=0.333\textwidth]{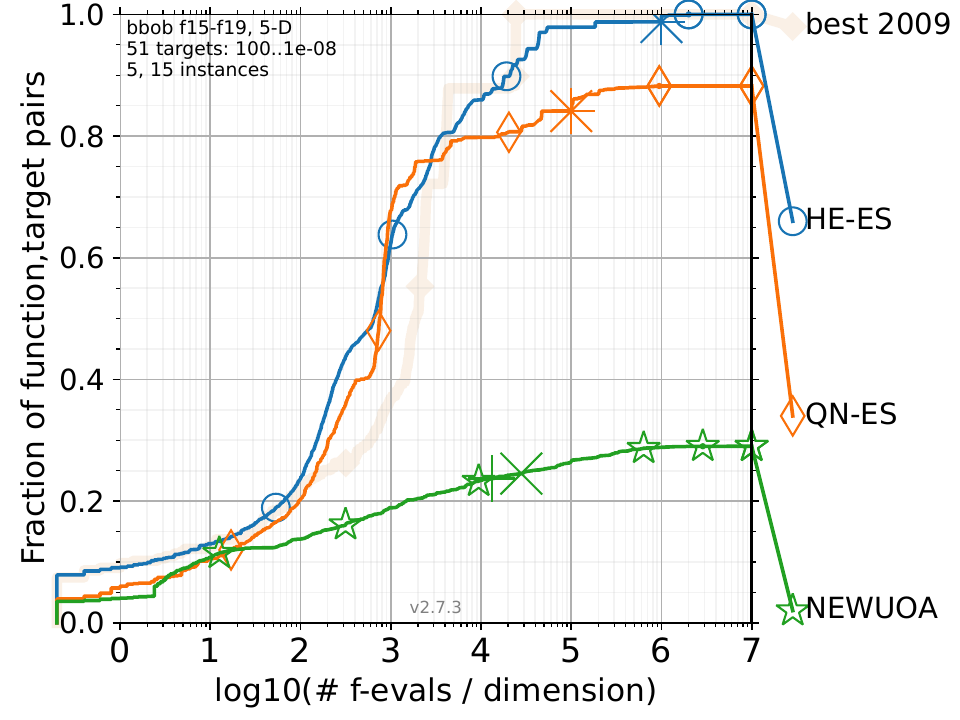}%
	\includegraphics[width=0.333\textwidth]{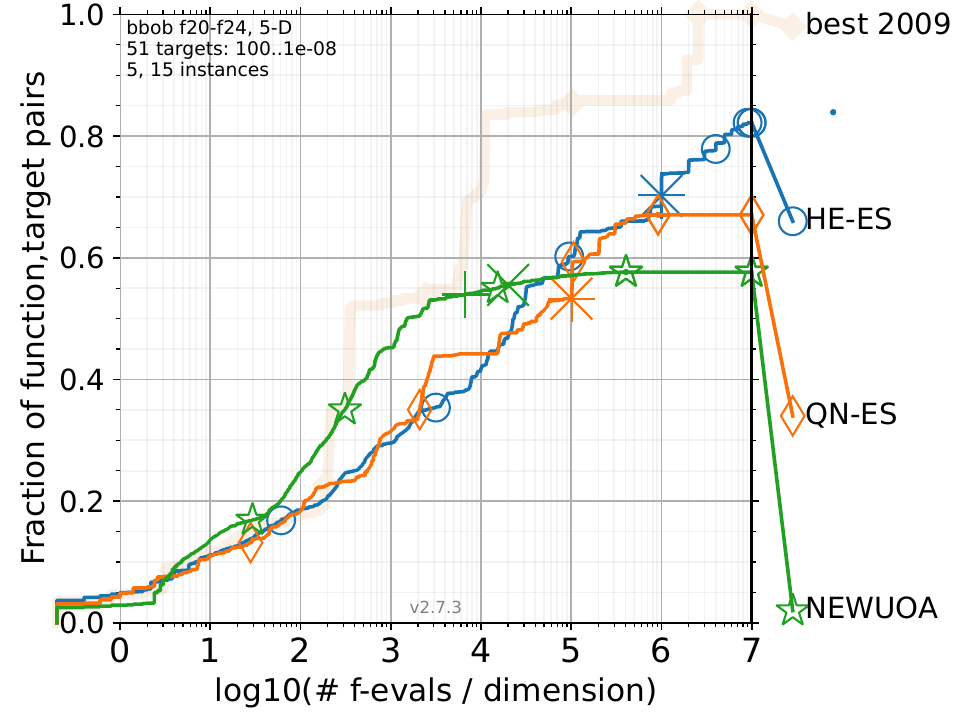}%
	\caption{
		\label{figure:BBOB-5D}
		COCO empirical cumulative distribution function (ECDF) plots by
		BBOB function group in problem dimension~5. The plots show the
		fraction of reached targets over the number of function
		evaluations divided by dimension (higher/more to the left is
		better). Note that beyond the large cross, the curves are
		extrapolated (based on ``virtual restarts'').
	}
\end{center}
\end{figure}

\begin{figure}
\begin{center}
	\includegraphics[width=0.333\textwidth]{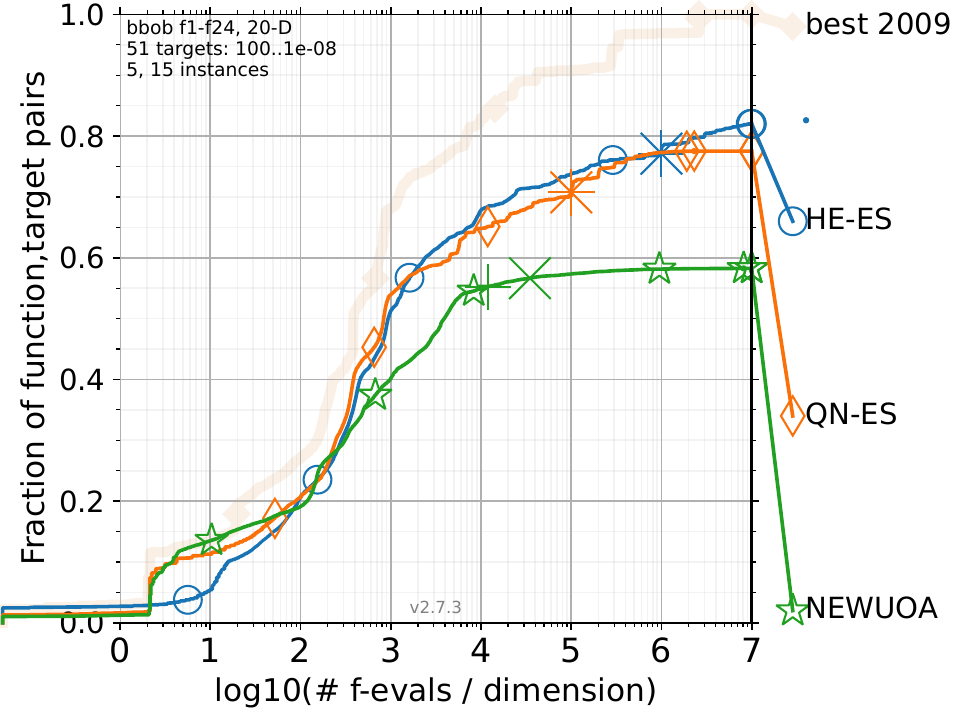}%
	\includegraphics[width=0.333\textwidth]{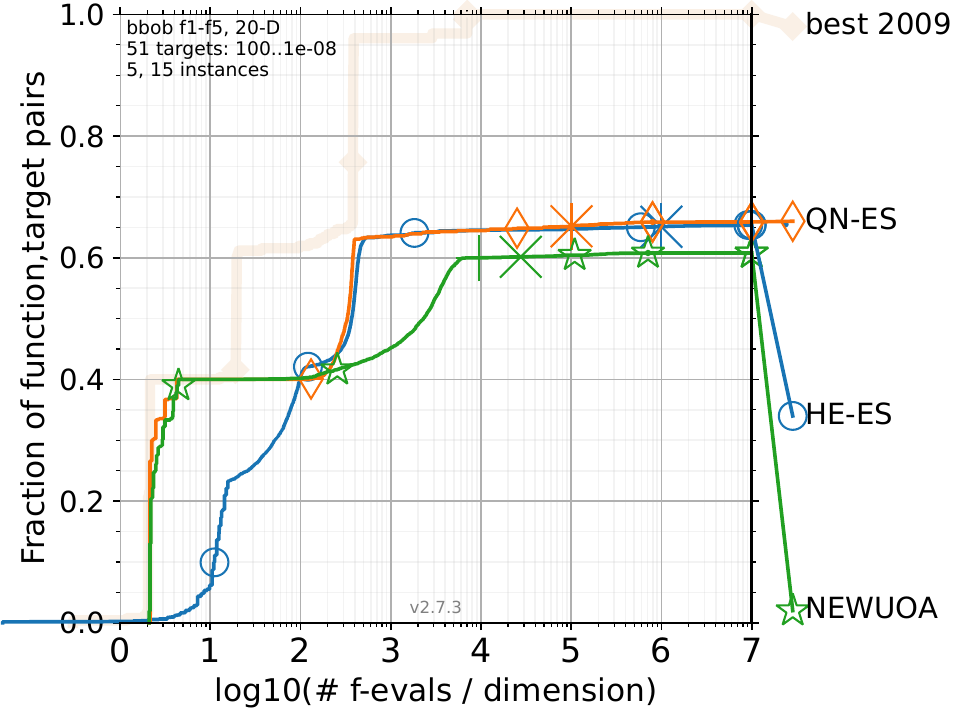}%
	\includegraphics[width=0.333\textwidth]{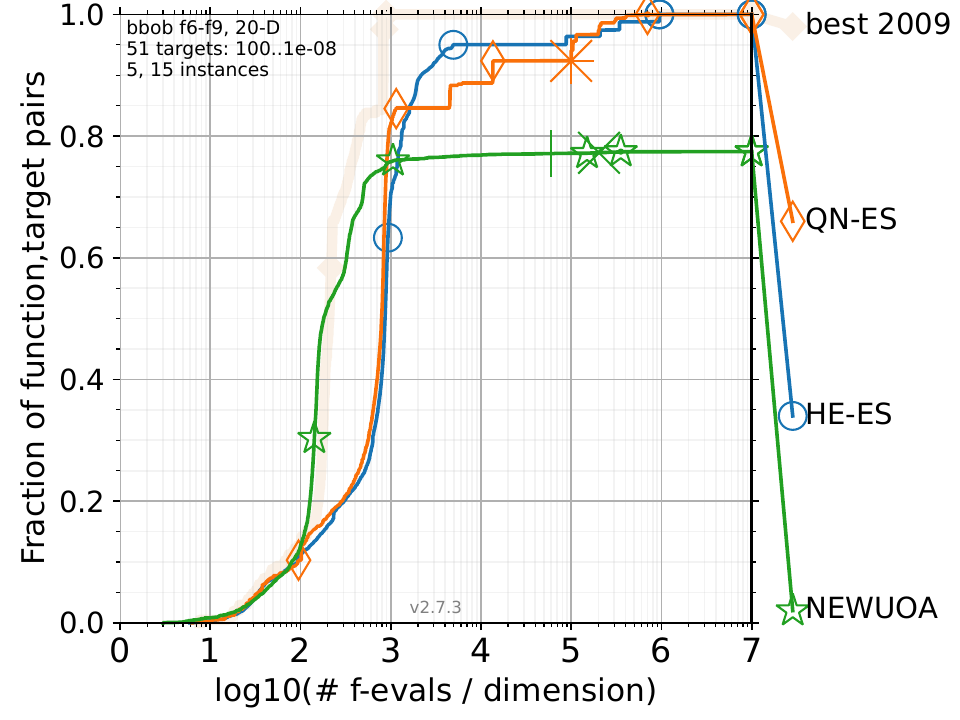}%
	\\
	\includegraphics[width=0.333\textwidth]{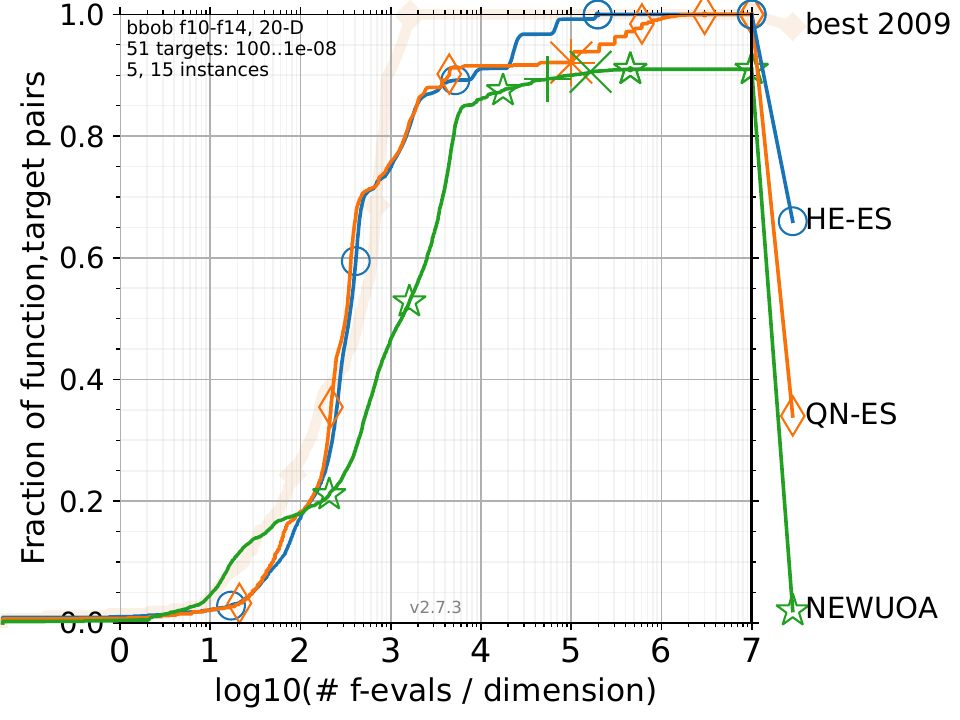}%
	\includegraphics[width=0.333\textwidth]{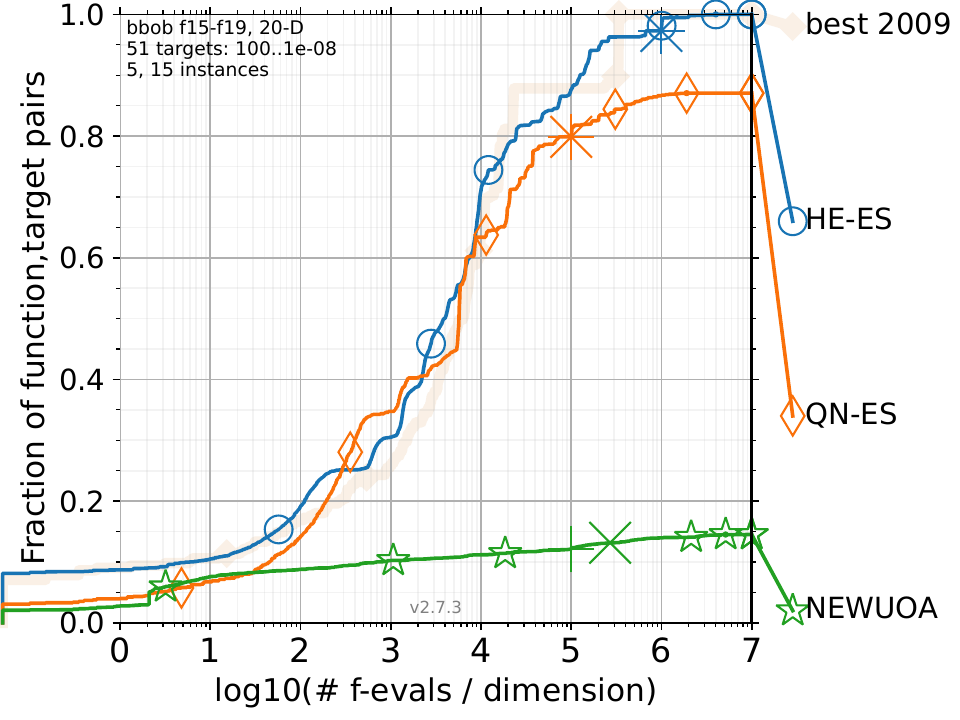}%
	\includegraphics[width=0.333\textwidth]{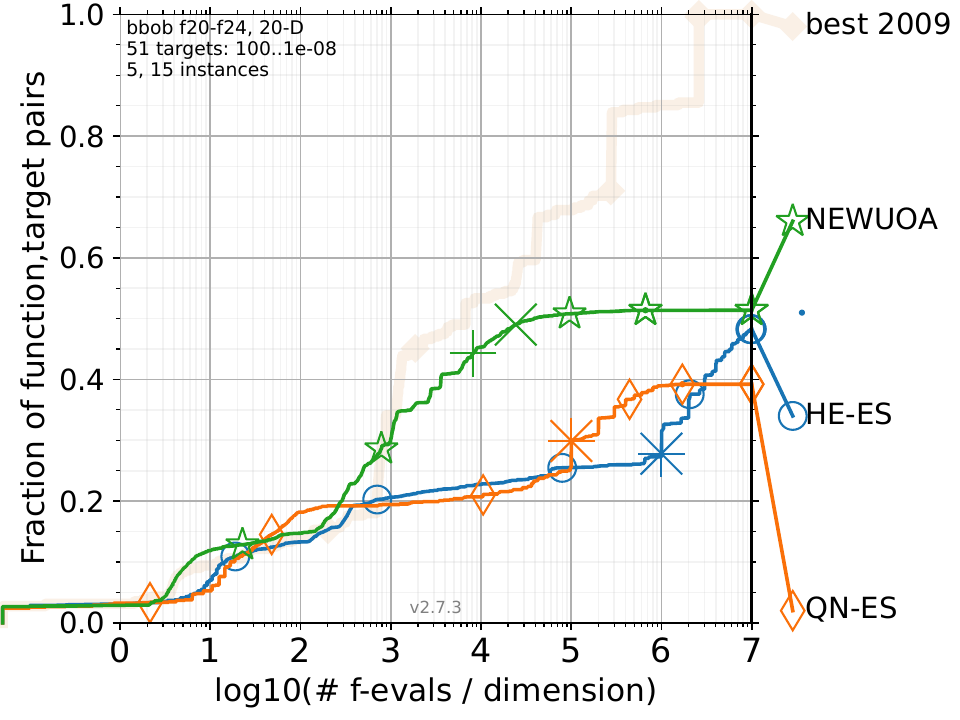}%
	\caption{
		\label{figure:BBOB-20D}
		COCO empirical cumulative distribution function (ECDF) plots by
		BBOB function group in problem dimension~20. The plots show the
		fraction of reached targets over the number of function
		evaluations divided by dimension (higher/more to the left is
		better). Note that beyond the large cross, the curves are
		extrapolated (based on ``virtual restarts'').
	}
\end{center}
\end{figure}

Results%
\footnote{
	In the spirit of open and reproducible research, the complete
	data set (including the other dimensions and performance data for each
	single function) is submitted to the COCO data repository.}
%% TODO: do it!
for dimensions 5 and 20 are presented in Figures \ref{figure:BBOB-5D}
and \ref{figure:BBOB-20D}. The top-left sub-figure aggregates
performance over all 24 problems. The remaining sub-figures refer to the
five function groups of separable problems (f1--f5), problems with mild
(f6--f9) and high (f10--f14) conditioning, as well as multi-modal
problems with (f15--f19) and without (f20--f24) global structure.

In accordance with the results presented above, the BBOB benchmark
indicates that QN-ES behaves more similarly to HE-ES than to NEWUOA. The
first sub-plot of each of the figures indicates that the evolution
strategies perform a little better overall on the (rather difficult)
BBOB benchmark.

Variants of many of the function listed in Table~\ref{table:benchmarks}
are included in BBOB, but they undergo slight non-linear
transformations. Although these transformations are sufficiently
moderate to preserve overall characteristics of the problems, they do
turn convex quadratic problems into (still unimodal) non-quadratic
problems. The effect on NEWUOA is quite dramatic, as can be observed by
its relatively weak performance on the group of ill-conditioned problems
f10--f14. In contrast, the two evolution strategies are nearly
unaffected. That also holds for QN-ES because its superlinear
convergence barely kicks in on the BBOB benchmark.

The performance on multi-modal problems seems to be governed by the
restart strategy. While HE-ES and QN-ES implement restarts with growing
population size, NEWUOA has no such means available. That turns out to
be a decisive deficit on problems f15--f19, but an advantage on problems
f20--f24. In accordance with \cite{hansen2010black} we assume that the
reason is that the very large number of restarts caused by short runs of
NEWOUA causes good performance on some of the problems. Analyzing the
performance of single functions (not shown), differences between HE-ES
and QN-ES can be found in a few cases. On f19, HE-ES performs much
better across all dimensions, while QN-ES outperforms HE-ES on f23. The
reasons for these differences remain opaque.

We do not notice significant dimension-dependent effects, with one
exception: QN-ES outperforms HE-ES on the separable function group in 5
dimensions, but not in 20 dimensions. A closer investigation reveals
that the difference stems entirely from the separable Rastrigin function
f3, which is solved by QN-ES up to dimension 10, while HE-ES and NEWUOA
start to degrade from dimension 5 onwards. The performance data does not
provide any hint at the underlying reason for that difference, and a
slight opposite effect shows up on the non-separable Rastrigin function
f15. That is unexpected, because all benchmarked optimizers are
invariant to orthogonal transformations of the search space. Neither
HE-ES nor QN-ES are designed to excel on the highly multi-modal
Rastrigin problem. Therefore the effect should not be overrated.

\section{Discussion and Conclusion}

% summary of the achievements
We have presented the (to the best of our knowledge) first evolution
strategy with superlinear convergence. It is a hybrid between the HE-ES
algorithm and a quasi-Newton DFO method. Its main algorithmic novelty is
to replace the well-established update of the distribution mean using
global weighted recombination with a quasi-Newton step.

% discussion of the results
Overall, our experimental results confirm the expectation that QN-ES
behaves mostly like an ES rather than a DFO method, but it offers a
tradeoff between robustness and efficiency that moves slightly into the
playing field of DFO algorithms like NEWUOA. Quasi-Newton steps require
a rather reliable quadratic model (or good estimates of first and second
derivatives), an assumption that is easily violated on hard problems.
Then we observe that QN-ES, although quite robust, falls slightly behind
HE-ES. However, it possesses the conceptually appealing property of
superlinear convergence, which brings significant benefits on smooth
problems.

% outlook
We have presented a novel algorithm, but we neither offer a mathematical
analysis nor did we explore the whole design space opened up by the new
concept. We therefore believe that our study opens more questions than
it answers. One line of future research is on the theory side, namely on
formally establishing superlinear convergence of QN-ES or a variant
thereof. A very practical line of research is algorithm design, aiming
to close the gap between ES and DFO algorithms further.

%\bibliographystyle{plain}
%\bibliography{bibliography}

\begin{thebibliography}{10}

\bibitem{akimoto2018drift}
Youhei Akimoto, Anne Auger, and Tobias Glasmachers.
\newblock Drift theory in continuous search spaces: expected hitting time of
  the {(1+1)-ES} with 1/5 success rule.
\newblock In {\em Proceedings of the Genetic and Evolutionary Computation
  Conference}, pages 801--808, 2018.

\bibitem{akimoto2022global}
Youhei Akimoto, Anne Auger, Tobias Glasmachers, and Daiki Morinaga.
\newblock Global linear convergence of evolution strategies on more than smooth
  strongly convex functions.
\newblock {\em SIAM Journal on Optimization}, 32(2):1402--1429, 2022.

\bibitem{auger2005restart}
Anne Auger and Nikolaus Hansen.
\newblock A restart {CMA} evolution strategy with increasing population size.
\newblock In {\em IEEE Congress on Evolutionary Computation}, volume~2, pages
  1769--1776, 2005.

\bibitem{beyer2001theory}
Hans-Georg Beyer.
\newblock {\em {The Theory of Evolution Strategies}}.
\newblock Springer Science \& Business Media, 2001.

\bibitem{beyer2014convergence}
Hans-Georg Beyer.
\newblock Convergence analysis of evolutionary algorithms that are based on the
  paradigm of information geometry.
\newblock {\em Evolutionary Computation}, 22(4):679--709, 2014.

\bibitem{beyer2012happycat}
Hans-Georg Beyer and Steffen Finck.
\newblock Happycat--a simple function class where well-known direct search
  algorithms do fail.
\newblock In {\em International conference on parallel problem solving from
  nature}, pages 367--376. Springer, 2012.

\bibitem{conn2009introduction}
Andrew~R Conn, Katya Scheinberg, and Luis~N Vicente.
\newblock {\em Introduction to Derivative-Free Optimization}.
\newblock SIAM, 2009.

\bibitem{gissler2024linear}
Armand Gissler.
\newblock {\em Linear Convergence of Evolution Strategies with Covariance
  Matrix Adaptation}.
\newblock PhD thesis, {\'E}cole Polytechnique, 2024.

\bibitem{glasmachers2020hessian}
Tobias Glasmachers and Oswin Krause.
\newblock {The Hessian Estimation Evolution Strategy}.
\newblock In {\em Parallel Problem Solving from Nature (PPSN XVII)}. Springer,
  2020.

\bibitem{glasmachers2022convergence}
Tobias Glasmachers and Oswin Krause.
\newblock {Convergence Analysis of the Hessian Estimation Evolution Strategy}.
\newblock {\em Evolutionary Computation Journal (ECJ)}, 30(1):27--50, 2022.

\bibitem{hansen2010comparing}
Nikolaus Hansen, Anne Auger, Raymond Ros, Steffen Finck, and Petr
  Po{\v{s}}{\'\i}k.
\newblock Comparing results of 31 algorithms from the black-box optimization
  benchmarking bbob-2009.
\newblock In {\em Proceedings of the 12th annual conference companion on
  Genetic and evolutionary computation}, pages 1689--1696, 2010.

\bibitem{hansen2001completely}
Nikolaus Hansen and Andreas Ostermeier.
\newblock Completely derandomized self-adaptation in evolution strategies.
\newblock {\em Evolutionary computation}, 9(2):159--195, 2001.

\bibitem{hansen2010black}
Nikolaus Hansen and Raymond Ros.
\newblock Black-box optimization benchmarking of {NEWUOA} compared to
  {BIPOP-CMA-ES}.
\newblock In {\em Proceedings of the 12th annual Conference Companion on
  Genetic and Evolutionary Computation}, pages 1519--1526, 2010.

\bibitem{jaegerskupper2007algorithmic}
Jens J{\"a}gersk{\"u}pper.
\newblock Algorithmic analysis of a basic evolutionary algorithm for continuous
  optimization.
\newblock {\em Theoretical Computer Science}, 379(3):329--347, 2007.

\bibitem{kern2004learning}
Stefan Kern, Sibylle~D M{\"u}ller, Nikolaus Hansen, Dirk B{\"u}che, Jiri
  Ocenasek, and Petros Koumoutsakos.
\newblock Learning probability distributions in continuous evolutionary
  algorithms--a comparative review.
\newblock {\em Natural Computing}, 3(1):77--112, 2004.

\bibitem{krause2019large}
Oswin Krause.
\newblock Large-scale noise-resilient evolution-strategies.
\newblock In {\em Proceedings of the Genetic and Evolutionary Computation
  Conference}, pages 682--690, 2019.

\bibitem{morinaga2019generalized}
Daiki Morinaga and Youhei Akimoto.
\newblock Generalized drift analysis in continuous domain: linear convergence
  of {(1+1)-ES} on strongly convex functions with lipschitz continuous
  gradients.
\newblock In {\em Proceedings of the 15th ACM/SIGEVO Conference on Foundations
  of Genetic Algorithms}, pages 13--24, 2019.

\bibitem{nocedal2006numerical}
Jorge Nocedal and Stephen~J Wright.
\newblock {\em Numerical Optimization (2nd Edition)}.
\newblock Springer, 2006.

\bibitem{ollivier2017information}
Yann Ollivier, Ludovic Arnold, Anne Auger, and Nikolaus Hansen.
\newblock Information-geometric optimization algorithms: A unifying picture via
  invariance principles.
\newblock {\em Journal of Machine Learning Research}, 18(18):1--65, 2017.

\bibitem{ostermeier1994step}
Andreas Ostermeier, Andreas Gawelczyk, and Nikolaus Hansen.
\newblock Step-size adaptation based on non-local use of selection information.
\newblock In {\em Parallel Problem Solving from Nature (PPSN)}, pages 189--198.
  Springer, 1994.

\bibitem{powell2006newuoa}
Michael~JD Powell.
\newblock The {NEWUOA} software for unconstrained optimization without
  derivatives.
\newblock In {\em Large-scale nonlinear optimization}, pages 255--297.
  Springer, 2006.

\bibitem{ragonneau2023pdfo}
T.~M. Ragonneau and Z.~Zhang.
\newblock {PDFO}: a cross-platform package for {Powell}'s derivative-free
  optimization solvers.
\newblock Technical Report arXiv:2302.13246, arXiv.org, 2022.

\bibitem{rechenberg1973evolutionsstrategie}
Ingo Rechenberg.
\newblock Evolutionsstrategie.
\newblock {\em Optimierung technischer Systeme nach Prinzipien der biologischen
  Evolution}, 1973.

\bibitem{teytaud2006general}
Olivier Teytaud and Sylvain Gelly.
\newblock General lower bounds for evolutionary algorithms.
\newblock In {\em Parallel Problem Solving from Nature-PPSN IX}, pages 21--31.
  Springer, 2006.

\bibitem{wierstra2014natural}
Daan Wierstra, Tom Schaul, Tobias Glasmachers, Yi~Sun, Jan Peters, and
  J{\"u}rgen Schmidhuber.
\newblock Natural evolution strategies.
\newblock {\em The Journal of Machine Learning Research}, 15(1):949--980, 2014.

\end{thebibliography}

\end{document}